\documentclass[letterpaper,onecolumn,12pt]{article}    
\usepackage{amssymb}
\usepackage{amsmath}
\usepackage{graphicx}  
\usepackage{epstopdf}
\usepackage{hyperref}
\usepackage{color}
\usepackage{bbold}
\usepackage{units}
\usepackage{a4wide}

\newtheorem{theorem}{Theorem}
\newtheorem{corollary}{Corollary}
\newtheorem{lemma}{Lemma}
\newtheorem{remark}{Remark}
\newtheorem{assumption}{Assumption}

\newtheorem{definition}{Definition}
\newtheorem{example}{Example}

\definecolor{changecolor}{rgb}{0.0, 0.0, 0.0}

\newcommand{\defeq}{\overset{\mathrm{def}}{=}}
\newcommand{\tr}{\intercal}

\newcommand{\qed}{$\diamond$}

\title{\textbf{Convex operator-theoretic methods in stochastic control}}

\author{Boris Houska}
\date{\small ShanghaiTech University}

\begin{document}

\addtolength{\textheight}{-2cm}

\maketitle


\begin{abstract}
This paper is about operator-theoretic methods for solving nonlinear stochastic optimal control problems to global optimality. These methods leverage on the convex duality between optimally controlled diffusion processes and Hamilton-Jacobi-Bellman (HJB) equations for nonlinear systems in an ergodic Hilbert-Sobolev space. In detail, a generalized Bakry-Emery condition is introduced under which one can establish the global exponential stabilizability of a large class of nonlinear systems. It is shown that this condition is sufficient to ensure the existence of solutions of the ergodic HJB for stochastic optimal control problems on infinite time horizons. Moreover, a novel dynamic programming recursion for bounded linear operators is introduced, which can be used to numerically solve HJB equations by a Galerkin projection.
\end{abstract}

\section{Introduction}
Bellman's principle of optimality~\cite{Bellman1957} and its dual counterpart, the Chapman-Kolmogorov equation~\cite{Chapman1928,Kolmogoroff1931} for controlled diffusion processes~\cite{Krylov2008}, are the fundament of stochastic control theory. They describe the semigroup structure of the evolution of optimally operated stochastic processes. Their associated partial differential equations (PDEs) for continuous-time models are the well-known HJB equation~\cite{Fleming1993} and the controlled Fokker-Planck-Kolmogorov (FPK) equation~\cite{Bogachev2015}.

\bigskip
\noindent
While parabolic FPKs and HJBs are related to stochastic optimal control problems on finite time horizons, their elliptic counterparts model the evolution of optimally operated diffusion processes on an infinite horizon. One distinguishes between \textit{ergodic optimal control} and \textit{infinite horizon optimal control}.

\begin{itemize}

\item \textit{Ergodic optimal control} is about the analysis of the limit behavior of optimally operated dynamic systems as well as the long-term average control performance. For instance, in the deterministic case, a control system might have optimal steady states or optimal periodic orbits, but it could as well be that it is optimal to stay on a homoclinic orbit or even to allow a chaotic limit behavior. Similarly, the optimal ergodic distribution of nonlinear stochastic control systems often has a highly complex multi-modal structure.

\item \textit{Infinite horizon optimal control} is about optimizing the system's transient behavior. Here, the hope is that this transient eventually converges to an optimal ergodic limit. For general nonlinear stochastic systems, however, such a convergence analysis often turns out to be difficult. This is because, in practice, one is often interested in optimizing the system's economic performance, which leads to indefinite stage cost functions. Unfortunately, it is often not even clear whether the integral over such indefinite costs has an expected value on an infinite time horizon.

\end{itemize}

\bigskip
\noindent
The goal of the present paper is to address the above challenges in ergodic- and infinite horizon optimal control by building upon recent results in elliptic diffusion operator theory and functional analysis, reviewed below.

\subsection{Literature Review}
The systematic mathematical study of stochastic systems goes back to the work of Kolmogorov~\cite{Kolmogoroff1931}, although related continuous-time electro-dynamic, molecular, quantum-theoretic, and thermal diffusions have been studied even earlier by Fokker~\cite{Fokker1914}, von Smoluchowski \cite{Smoluchowski1915}, Planck~\cite{Planck1917}, and Chapman~\cite{Chapman1928}. The relevance of such diffusions in the context of analyzing the ergodic behavior of dynamical systems has been observed by Has'minski\v{i}~\cite{Hasminskii1960} in 1960, who found an elegant way of using Harnack's inequality in combination with a weak Lyapunov function in order to analyze the existence of steady states of FPK equations. And, in fact, it was only six years later that Wonham realized the great potential of Has'minski\v{i}'s work in the context of ergodic control theory~\cite{Wonham1966}. By now, these historical developments are well-understood and have been extended by many authors, as reviewed in great detail in~\cite{Bogachev2015} focusing on FPK equations and in~\cite{Krylov2008} focusing on controlled diffusion processes.

\bigskip
\noindent
In contrast to linear FPKs, HJBs are nonlinear PDEs, which have been introduced by Bellman~\cite{Bellman1957}. They model the backward-in-time evolution of the value function of either deterministic or stochastic optimal control problems. Here, the state of the HJB equation can---at least formally---be interpreted as the co-state of the FPK-constraint of certain optimal diffusion control problems. As such, HJBs and FPKs are closely related via the concept of duality in convex optimization, as originally pointed out in~\cite{Bismut1973} and~\cite{Fleming1989}. HJBs for optimal control of deterministic systems are typically analyzed by introducing an artificial stochastic perturbation and then passing to the so-called viscosity limit by making this perturbation infinitesimally small. The corresponding viscosity methods go back to Crandall, Ishii and Lions~\cite{Crandall1992}. They can be considered as the main tool for analyzing parabolic HJBs on finite time horizons. The textbooks~\cite{Bardi1997} and \cite{Fleming1993} offer a relatively complete overview of the associated literature on HJBs between 1957 and 1997.

\bigskip
\noindent
Regarding the existence and properties of solutions to elliptic HJBs for ergodic and infinite-horizon optimal control, however, surprisingly little is known. In~\cite{Bardi1997} a whole chapter is dedicated to infinite horizon optimal control problems with exponentially discounted stage costs. The analysis methods in this chapter are, however, inapplicable to problems with a general indefinite stage cost. In the context of ergodic optimal control, this problem has partly been resolved by Arisawa~\cite{Arisawa1998b,Arisawa1998}, who has proposed methods to analyze the limit of vanishing discounts, which, under certain assumptions, turns out to be a successful path towards analyzing ergodic HJBs. Moreover, Barles~\cite{Barles2000} as well as Namah and Roquejoffre~\cite{Namah1999} have suggested methods for analyzing the long term behavior of HJBs with periodic Hamiltonians. These methods have in common that they focus on special classes of stochastic systems under periodicity assumptions. And, finally, yet another option for analyzing infinite horizon optimal control problems is to work with turnpikes~\cite{Askovic2022,Trelat2015}, which can be used to analyze the long-term local behavior of systems near optimal steady-states or periodic orbits.

\bigskip
\noindent
In contrast to ergodic HJBs, much more is known about the limit behavior of FPKs. This is not only due the above mentioned work by Has'minski\v{i}, but also due to relatively recent breakthroughs in functional analysis that can be traced back to the seminal work of Bakry and Emery~\cite{Bakry1985} on \textit{carre-du-champ} operators, whose far-reaching consequences have only been understood to full extent in the last two decades. Based on earlier work by Gross on hypercontractivity and logarithmic Sobolev inequalities~\cite{Gross1975}, Arnold, Markowich, Toscani, and Unterreiter~\cite{Arnold2001} were among the first to come up with a comprehensive and general analysis of the convergence of FPKs to ergodic limit states by featuring Bakry-Emery conditions. Even more practical conditions for deriving logarithmic Sobolev- and general entropy-type inequalities for diffusion semigroups and FPKs have appeared in an article by Bolley and Gentil~\cite{Bolley2010}. The most recent developments regarding the trend to equilibrium of FPKs in the degenerate case are summarized in~\cite{Barbu2021,Barbu2023}.

\bigskip
\noindent
Last but not least, concerning recent developments in PDE-constrained optimization~\cite{Hinze2009}, optimal control of FPKs has received significant attention during the last five years~\cite{Aronna2021,Breiten2018,Breiten2018b}. These articles are, however, more concerned with the well-posedness of FPK-constrained optimization problems, typically, with least-squares or other objectives. These particular classes of FPK-constrained optimization problems are not directly related to HJBs. The only exception is the work of Bobkov, Gentil, and Ledoux~\cite{Bobkov2001}, who are the first to attempt an analysis of Hamilton-Jacobi equations based on logarithmic Sobolev inequalities---even if these authors focus only on a very special class of HJBs that can be solved explicitly via the Hopf-Lax formula.

\subsection{Contributions}
This paper presents novel results about the existence and uniqueness of solutions to ergodic and infinite horizon stochastic optimal control problems. They are summarized, respectively, in Theorem~\ref{thm::ergodic} and Theorem~\ref{thm::HJB}.

\bigskip
\noindent
In detail, Theorem~\ref{thm::ergodic} is about the existence of ergodic limits of optimally controlled nonlinear stochastic McKean-Vlasov systems on a general convex domain with reflecting boundary, assuming that a weak Has'minski\v{i}-Lyapunov function is available. Its proof is based on an energy dissipation argument for FPKs with mixed boundary conditions that models the evolution of the probability density function of its underlying nonlinear stochastic control system. The theorem is applicable to a general class of ergodic optimal control problems with indefinite stage cost and potentially multi-modal optimal limit distributions.

\bigskip
\noindent
Similarly, Theorem~\ref{thm::HJB} is about the existence of weak solutions to ergodic HJBs with Neumann boundary condition in a suitably weighted Hilbert-Sobolev space and their relation to infinite horizon stochastic optimal control. It is based on a generalized Bakry-Emery condition for optimally controlled FPKs under which global exponential stabilizability can be established. This condition turns out to be sufficient for ensuring that the value function of the stochastic infinite horizon optimal control problem is a bounded linear operator, whose unique Riesz' representation can be used to construct a weak solution of the ergodic HJB.

\bigskip
\noindent
A particularly relevant aspect of Theorem~\ref{thm::HJB} is that it succeeds in analyzing a class of nonlinear \mbox{PDEs---in} this case HJBs---by relying on convex functional analysis methods that have originally been designed for analyzing linear PDEs. This implies in particular that mature numerical methods for linear parabolic PDEs, such as finite element Galerkin methods, can be used to approximately solve finite- and infinite horizon stochastic optimal control problems via a discretized operator-theoretic dynamic programming recursion. Because such methods for discretizing linear parabolic PDEs are, in general, more accurate, easier to analyze, and easier to implement than methods for discretizing nonlinear PDEs, such as HJBs, operator-theoretic dynamic programming methods offer clear advantages over classical dynamic programming methods for optimal control.

\subsection{Overview}
In order to put the above contributions into context, this paper is structured as follows.

\begin{itemize}
\addtolength{\itemsep}{2pt}

\item Section~\ref{eq::scs} reviews the relations between controlled McKean-Vlasov systems and FPKs in the presence of reflecting boundary conditions.

\item Section~\ref{sec::ergodic} explains how to use energy functionals to establish dissipativity results for optimally controlled FPKs, leading to the first main result on ergodic optimal control that is summarized in Theorem~\ref{thm::ergodic}.

\item Section~\ref{sec::HJB} proposes a novel approach to analyzing parabolic HJBs with Neumann boundary conditions in a suitably weighted Hilbert-Sobolev space.

\item Section~\ref{sec::BEC} introduces Lyapunov diffusion operators in order to  analyze the stabilizability of controlled diffusions via a generalized Bakry-Emery condition.

\item Section~\ref{sec::IHOCP} introduces ergodic HJBs and presents the main result on infinite horizon optimal control that is summarized in Theorem~\ref{thm::HJB}.

\item Section~\ref{sec::algorithm} presents a finite element Galerkin projection based operator-theoretic dynamic programming method for solving parabolic and ergodic HJBs. A case study illustrates its performance for an infinite horizon stochastic optimal control problem with indefinite stage costs and non-trivial ergodic behavior.

\end{itemize}

\bigskip
\noindent
And, finally, Section~\ref{sec::conclusions} concludes the paper.

\subsection{Notation}
Let $\Omega \subseteq \mathbb R^{\mathsf{n}}$ denote an open set, $\overline \Omega = \mathrm{cl}(\Omega)$ its closure and $\partial \Omega = \overline \Omega \setminus \Omega$ its boundary. The sets
\[
L_\mathrm{loc}^k(\Omega), \ L^k(\Omega), \ C_0^\infty(\Omega), \ C^k \left( \overline \Omega \right), \ \text{and} \ \ H^1(\Omega)
\]
denote---as usual---the set of functions that are, respectively, locally $L^k$-integrable, $L^k$-integrable, smooth with compact support, $k$-times continuously differentiable with continuous extensions to the boundary, and $L^2$-integrable with $L^2$-integrable weak derivatives on $\Omega$. The inner products
\begin{eqnarray}
( \varphi, \psi )_{L^2(\Omega)} & \ \defeq \ & \int_{\Omega} \varphi \psi \, \mathrm{d}x  \notag \\[0.16cm]
\text{and} \qquad ( \varphi, \psi )_{H^1(\Omega)} & \ \defeq \ & \int_{\Omega} \left( \varphi \psi + \nabla \varphi^\tr \nabla \psi \right) \, \mathrm{d}x \notag
\end{eqnarray}
are defined such that $L^2(\Omega)$ and $H^1(\Omega)$ are Hilbert spaces. Here, $\nabla$ denotes the gradient operator. Additionally, if $\omega \in L_\mathrm{loc}^1(\Omega)$ with $\omega > 0$ is a strictly positive weight function, the shorthands
\begin{eqnarray}
L^2(\omega) & \ \defeq \ & \left\{ \ \varphi \in L_\mathrm{loc}^1(\Omega) \ \middle| \ \varphi \cdot \sqrt{\omega} \ \in \ L^2(\Omega) \ \right\} \notag \\[0.16cm]
\text{and} \quad H^1(\omega) & \ \defeq \ & \left\{ \ \varphi \in L_\mathrm{loc}^1(\Omega) \ \middle| \ \varphi \cdot \sqrt{\omega} \ \in \ H^1(\Omega) \ \right\} \notag 
\end{eqnarray}
are used. The associated inner products are then, of course, weighted with $\omega$, too. For given time intervals $[0,T] \subseteq \mathbb R$ and any Banach space $(X,\|\cdot\|_X)$, we use the symbol $L^k(0,T;X)$ to denote the set of strongly measurable functions $\varphi: [0,T] \to X$ with
\[
\| \varphi \|_{L^k(0,T;X)} \ \defeq \ \left( \int_0^T \| \varphi(t) \|_X^k \, \mathrm{d}t \right)^{1/k} \ < \ \infty ,
\]
such that $L^k(0,T;X)$ is a Bochner space. Weak derivatives with respect to time are indicated by a ``dot''; and $H^1(\omega)^*$ denotes the dual of $H^1(\omega)$. The space
\[
\mathcal P(\omega) \ \defeq \ \{ \ \varphi \in L^2(0,T;H^1(\omega)) \mid \dot \varphi \in L^2(0,T;H^1(\omega)^*) \ \}
\]
is then again a Hilbert space with respect to its natural inner product~\cite[Thm.~1.31]{Hinze2009}. Finally, the natural duality pairing between a Hilbert space and its dual is denoted by $\langle \cdot, \cdot \rangle$ assuming that it is clear from the context which Hilbert space is meant.

\section{Stochastic Control Systems}
\label{eq::scs}

This section is about the equivalence of closed-loop controlled stochastic nonlinear systems in a finite dimensional space and a special class of open-loop controlled deterministic linear systems in an infinite dimensional space, namely, controlled FPKs.

\subsection{Nonlinear Control Systems}
Deterministic control-affine systems have the form
\begin{eqnarray}
\label{eq::system}
\dot x = f(x) + G(x) u \qquad \text{on} \quad \mathbb T \ \defeq \ (0,T) \; .
\end{eqnarray}
Here, $\mathbb T \subseteq \mathbb R$ denotes an open time horizon. It is finite if $T < \infty$ and infinite if $T = \infty$. The right-hand functions $f \in C^1(\overline \Omega)^{n_x}$ and $G \in C^1(\overline \Omega)^{n_u}$ are assumed to be defined on the closure of an open domain $\Omega \subseteq \mathbb R^{n_x}$. Moreover, $x(t)$ denotes the state and $u(t)$ the control at time $t \in \mathbb T$.

\bigskip
\noindent
The goal is to design a feedback $\mu \in L_\mathrm{loc}^1(0,T;L_\mathrm{loc}^1(\Omega))$ that can be used to control~\eqref{eq::system} in closed-loop mode,
\begin{align}
\label{eq::ClosedLoop}
\dot x(t) = F[\mu](t,x(t)) \quad \text{for} \ \ t \in \mathbb T,
\end{align}
with $F[\mu](t,x) \ \defeq \ f(x) + G(x) \mu(t,x)$. Because traditional existence theorems for ordinary differential equations are inapplicable in the context of such general feedback functions $\mu$, solutions of~\eqref{eq::ClosedLoop} are defined in the viscosity sense, by introducing a small random process noise.

\subsection{Stochastic Differential Equations}
Let us temporarily assume that $\Omega = \mathbb R^{n_x}$. In this case, we can add a process noise to the right-hand of~\eqref{eq::ClosedLoop}, arriving at a stochastic differential equation (SDE) of the form
\begin{eqnarray}
\label{eq::SDE}
\mathrm{d} X = F[\mu](t,X) \, \mathrm{d}t + \sqrt{2\epsilon} \, \mathrm{d} W_t \qquad \text{on} \ \ \mathbb T,
\end{eqnarray}
where $W_t$, $t \geq 0$, denotes a Wiener process, $\epsilon > 0$ a diffusion constant, and $x_0 \in \mathbb R^{n_x}$ an initial state. The state of this SDE at time $t \in \mathbb T$ is denoted by $X(t)$. It is a random variable, whose probability density function is denoted by $\rho(t) \in L^1(\Omega)$. This means that for any given Borel set $S \subseteq \Omega$, the probability of the event $X(t) \in S$ is given by
\[
\mathrm{Pr} \left( \, X(t) \in S \, \right) \ = \ \int_{S} \rho(t) \, \mathrm{d}x \; .
\]
Conditions under which~\eqref{eq::SDE} admits a martingale solution $X$ that has such an $L^1$-integrable probability density function $\rho$ are discussed below.

\begin{remark}
The assumption that $\epsilon$ is a scalar constant is introduced for simplicity presentation. The results below can be generalized to the case that the constant scalar factor $\sqrt{2 \epsilon}$ in~\eqref{eq::SDE} is replaced by a potentially state-dependent invertible matrix.
\end{remark}

\subsection{Reflections at the Boundary}
The states of many nonlinear systems have no physical interpretation outside of their natural domain $\Omega$. For instance the concentration of the reactants in a controlled chemical process are usually required to be non-negative implying that $\Omega \neq \mathbb R^{n_x}$. In such a case, the stochastic process noise model in~\eqref{eq::SDE} makes no sense, as the Wiener process $W_t$ is unbounded. Instead, at least for the case that $\Omega$ is convex, a physically more meaningful stochastic model is given by
\begin{eqnarray}
\label{eq::SDE2}
\begin{array}{rcl}
\sqrt{2 \epsilon} \, \mathrm{d}W_t &\in& \mathrm{d} X - F[\mu](t,X) \mathrm{d}t + N_{\overline \Omega}(X) \, \mathrm{d}t \quad \text{on} \ \ \mathbb T . \ \ \
\end{array}
\end{eqnarray}
This is a McKean-Vlasov SDE with reflection~\cite{Ciotir2021,Sznitman1984}, where the additional drift term,
\[
N_{\overline \Omega}(x) \ \defeq \ \{ z \in \mathbb R^{n_x} \mid \forall y \in \overline \Omega, \ z^\tr(x-y) \geq 0 \},
\]
acts as reflector. We have $N_{\overline \Omega}(x) = \{ 0 \}$ for $x \in \Omega$, which means that~\eqref{eq::SDE} and~\eqref{eq::SDE2} coincide on $\Omega$. For $x \in \partial \Omega$, however, $N_{\overline \Omega}(x)$ is equal to the normal cone of $\overline \Omega$ at $x$. Consequently, whenever~\eqref{eq::SDE2} admits a solution $X$ that has an $L^1$-integrable probability density $\rho$, we expect that
\begin{align}
\label{eq::prob1}
\mathrm{Pr}( \, X(t) \in \Omega \, ) \ = \ \int_{\Omega} \rho(t) \, \mathrm{d}x \ = \ 1 \; ,
\end{align}
which is the case under very mild assumptions that are, for instance, discussed in~\cite{Ciotir2021}, but sufficient conditions under which this is the case are also discussed below. For the special case that $\Omega = \mathbb R^{n_x}$,~\eqref{eq::SDE2} reduces to~\eqref{eq::SDE}.

\subsection{Fokker-Planck-Kolmogorov Equations}

Let $\mathcal A$ denote the $\epsilon$-diffusion operator of $f$ and $\mathcal B$ the linear transport operator of $G$, formally defined by
\[
\mathcal A \rho = -\mathrm{div}\left( f \rho \right) + \epsilon \, \Delta \rho \qquad \text{and} \qquad \mathcal B \nu = -\mathrm{div}(G \nu) \; ,
\]
where ``$\mathrm{div}$'' denotes the divergence and ``$\Delta$'' the Laplace operator. The pair $(\mathcal A,\mathcal B)$ is associated with an infinite-dimensional linear control system, called the controlled FPK equation. In the presence of reflecting boundary conditions, it has the form
\begin{eqnarray}
\label{eq::FokkerPlanck}
\begin{array}{rcll}
\dot \rho &=& \mathcal A \rho + \mathcal B \nu & \text{on} \ \ \Omega_{\mathbb T} \\[0.16cm]
\rho_0 &=& \rho(0) & \text{on} \ \ \Omega \\[0.16cm]
0 &=& \left. \left( \epsilon \nabla \rho - f \rho - G \nu \right)^\tr n \right|_{\partial \Omega} \quad & \text{on} \ \ \Sigma_{\mathbb T}. \\
\end{array}
\end{eqnarray}
Here, $\rho_0 \in L^1(\Omega)$ denotes a given initial probability distribution on $\Omega$, while \[
\Omega_{\mathbb T} \ \defeq \ \mathbb T \times \Omega \qquad \text{and} \qquad
\Sigma_\mathbb T \ \defeq \ \mathbb T \times \partial \Omega
\]
denote, respectively, the FPK's open domain and its spatial boundary. Moreover, \mbox{$n: \partial \Omega \to \mathbb R^{n_x}$} denotes an outer normal, such that
\[
\forall x \in \partial \Omega, \quad n(x) \in N_{\overline \Omega}(x) \quad \text{and} \quad \| n(x)\|_2 = 1 \; . 
\]
The open-loop linear control system~\eqref{eq::FokkerPlanck} is closely related to the closed-loop nonlinear control system~\eqref{eq::SDE2} after identifying $\rho$ with the probability distribution of the random process $X$ and after substituting
\begin{eqnarray}
\label{eq::murho}
\nu \ = \ \mu \rho \qquad \text{on} \ \ \Omega_{\mathbb T} \; .
\end{eqnarray}
In order to elaborate on this relation between~\eqref{eq::SDE2},~\eqref{eq::FokkerPlanck} and~\eqref{eq::murho}, a couple of definitions are needed. For this aim, the shorthand
\[
\mathcal U(\omega) \ \defeq \ (L^2(0,T;L^2(\omega)))^{n_u}
\]
is introduced, where $\omega \in L_\mathrm{loc}^1(\Omega)$, $\omega > 0$, is a strictly positive function. Moreover, the linear operators $\mathcal A$ and $\mathcal B$ are associated with the bilinear forms
\begin{eqnarray}
a(\rho,V) & \ \defeq \ & \int_{\Omega} \rho f^\tr \nabla V - \epsilon \nabla \rho^\tr \nabla V \, \mathrm{d}x \notag \\[0.2cm]
\text{and} \qquad b(\nu,V) & \ \defeq \ & \sum_{i=1}^{n_u} \left\langle \nu_i, \nabla V^\tr G_i \right\rangle, \notag
\end{eqnarray}
which are defined for all $\rho \in \mathcal P(\omega^{-1})$ and all $V \in \mathcal P(\omega)$.

\begin{definition}
\label{def::weak}
Let $\omega \in L^\infty(\Omega)$ with $\omega > 0$ on $\Omega$ and $\nu \in \mathcal U(\omega^{-1})$ be given. A function $\rho \in \mathcal P(\omega^{-1})$ with given initial value $\rho(0) = \rho_0 \in L^2(\omega^{-1})$ is called a weak solution of~\eqref{eq::FokkerPlanck} if
\begin{eqnarray}
\label{eq::weak1}
\int_0^T \left\langle \dot \rho, V \right\rangle \, \mathrm{d}t & \ = \ & \int_0^T a(\rho,V)  + b(\nu,V) \, \mathrm{d}t
\end{eqnarray}
for all test functions $V \in \mathcal P(\omega)$.
\end{definition}

\bigskip
\noindent
For the special case that $\Omega$ is convex and bounded and $T < \infty$, we may set the weight function to $\omega = 1$. In this case, the existence of a unique weak solution of~\eqref{eq::FokkerPlanck} for arbitrary $\nu \in \mathcal U(1)$ follows from a standard result about parabolic PDEs, see~\cite[Thm.~1.33]{Hinze2009}, recalling that $f$ and $G$ are assumed to be continuously differentiable on $\overline \Omega$ and $\epsilon > 0$. For unbounded domains $\Omega_{\mathbb T}$, however, additional assumptions on $f$, $G$, and $\omega$ are needed in order to ensure the existence of weak solutions, as discussed a bit further below, in Section~\ref{sec::ergodic}. 

\subsection{Nonlinear Control via Linear Systems Theory}
From a control theoretic perspective, it is interesting to point out that for the purpose of analyzing the finite-dimensional closed-loop controlled nonlinear system $(f,G)$, as defined by~\eqref{eq::SDE2}, it sufficient to analyze the infinite dimensional open-loop controlled linear system $(\mathcal A, \mathcal B)$, as defined in~\eqref{eq::FokkerPlanck}. This is because the existence of a unique weak solution $\rho \in \mathcal P(\omega^{-1})$ of~\eqref{eq::FokkerPlanck} for $\omega \in L^\infty(\Omega)$ and an initial probability density function $0 \leq \rho_0 \in L^2(\omega^{-1}) \cap L^1(\Omega)$ with $\int_{\Omega} \rho_0 \mathrm{d}x = 1$ of the initial random variable $X(0)$ is sufficient to ensure that~\eqref{eq::SDE2} has a martingale solution $X$ with probability density function $\rho$ that satisfies~\eqref{eq::prob1}. This is a standard result from the theory of diffusion processes~\cite{Tresivan2016}, which has recently been extended for general FPKs with reflecting boundary conditions~\cite{Barbu2021}; see also~\cite{Ciotir2021}. Of course, here we have to substitute~\eqref{eq::murho}, but this substitution is justified, too. Namely, it follows from the minimum principle for parabolic FPKs~\cite{Bogachev2015} that if $\rho \in \mathcal P(\omega^{-1})$, is a weak solution of~\eqref{eq::FokkerPlanck} for an initial probability density $\rho_0 \geq 0$, then
\begin{align}
\label{eq::minPrinciple}
\rho > 0 \quad \text{on} \ \ \ \Omega_{\mathbb T} \; ,
\end{align}
recalling that $\Omega_{\mathbb T}$ is an open set, $0 \notin \mathbb T$, and $\epsilon > 0$. Thus, the closed-loop control law, $\mu = \nu/\rho$, can be recovered from the infinite-dimensional open-loop input $\nu$ and the strictly positive function $\rho$.

\begin{remark}
There are many ways to reformulate finite dimensional nonlinear control systems as equivalent infinite dimensional linear control systems. This includes linear embedding methods based on occupation measures~\cite{Lasserre2008,Vinter1993} as well as lifting methods based on first order differential operators such as Liouville-Koopman or Perron-Frobenius operators~\cite{Korda2018,Villanueva2021}. In contrast to these first order transport operators and their adjoints, the Kolmogorov forward operator $\mathcal A$ is, however, an elliptic second order operator. Such second order operators have many advantages regarding both their mathematical properties and their numerical discretization~\cite{Hinze2009}.
\end{remark}

\section{Ergodic Optimal Control}
\label{sec::ergodic}
The idea to analyze optimal control of diffusion processes with reflection can be traced back to the work of Puterman~\cite{Puterman1977}, who originally introduced such optimal control problems on a finite time horizon $\mathbb T$ and on a bounded domain $\Omega$. This section presents methods for analyzing such optimal control problems on potentially unbounded domains, with a focus on ergodic limits.

\subsection{Problem Formulation}
Let $\mathbb U \subseteq \mathbb R^{n_u}$ denote a control constraint set, let
\[
\mathfrak U \ \defeq \ \{ \mu \in L_\mathrm{loc}^1(0,T;L_\mathrm{loc}^1(\Omega)) \mid \mu \in \mathbb U \ \ \text{on} \ \ \Omega_{\mathbb T} \}
\]
denote the associated set of feasible control laws, and let $\ell: \Omega \times \mathbb U \to \mathbb R$ be a cost function. The current section is about ergodic optimal control problems of the form
\begin{eqnarray}
\label{eq::EOC}
J^\star \ = \ \lim_{T \to \infty} \min_{\mu \in \mathfrak U} \  \mathbb{E} \left\{ \frac{1}{T} \int_0^T \ell(X(t), \mu(X(t)) ) \, \mathrm{d}t \right\} .
\end{eqnarray}
Here, $X$ denotes a martingale solution of the SDE~\eqref{eq::SDE2} that depends implicitly on the optimization variable $\mu$.

\begin{assumption}
\label{ass::blanket}
Throughout this paper, we assume that
\begin{enumerate}

\item the set $\Omega \subseteq \mathbb R^{n_x}$ is convex, open, and non-empty;

\item the system satisfies $(f,G) \in C^1(\overline \Omega)^{n_x \times (1+n_u)}$;

\item the set $\mathbb U \subseteq \mathbb R^{n_u}$ is convex, compact, $\mathbb U \neq \varnothing$; and

\item the stage cost $\ell \in C^2(\overline \Omega \times \mathbb U)$ is bounded from below, radially unbounded in $x$,
\[
\| x\| \to \infty \quad \Longrightarrow \quad \ell(x,u) \to \infty,
\]
and strongly convex in $u$, $\nabla_u^2 \ell \succ 0$ on $\Omega \times \mathbb U$.
\end{enumerate}

\end{assumption}
Assumption~\ref{ass::blanket} neither excludes that $\Omega$ is unbounded nor that $\Omega = \mathbb R^{n_x}$. On such unbounded domains, however, a weak Has'minski\v{i}-Lyapunov function is needed.

\begin{assumption}
\label{ass::Q}
There exists a function $Q \in C^2(\overline \Omega)$ and positive constants $\gamma_1,\gamma_2,\gamma_3,\gamma_4 \in (0,\infty)$ such that
\begin{enumerate}
\addtolength{\itemsep}{2pt}

\item the Hessian of $Q$ is positive definite and uniformly bounded, such that
\begin{eqnarray}
\label{eq::HessianBound}
\gamma_1 I \ \preceq \ \nabla^2 Q \ \preceq \ \gamma_2 I \qquad \text{on} \ \ \Omega ,
\end{eqnarray}

\item the function $Q$ satisfies $\nabla Q^\tr n \geq 0$ on $\partial \Omega$, and

\item  $Q$ satisfies the weak Lyapunov condition
\begin{eqnarray}
\label{eq::RelaxedLyapunov}
F[\mu]^\tr \nabla Q & \ \leq \ & \gamma_3 - \gamma_4 \ell[\mu] \qquad \text{on} \ \ \Omega \\[0.16cm]
\label{eq::ellmu}
\text{with} \qquad \quad \ell[\mu](x) &  \ \defeq \  & \ell(x,\mu(x))
\end{eqnarray}
for at least one $\mu \in C^1(\Omega)$ with $\mu \in \mathrm{int}(\mathbb U)$ on $\Omega$.

\end{enumerate}
\end{assumption}

\bigskip
\noindent
Assumption~\ref{ass::Q} can be interpreted as a weak stabilizability condition of the system $(f,G)$ relative to its cost $\ell$.

\begin{example}
If Assumption~\ref{ass::blanket} holds and if the domain $\Omega = \{ x \mid \| x \| < 1 \}$ is the open unit ball, Assumption~\ref{ass::Q} is always satisfied, because Condition~\eqref{eq::RelaxedLyapunov} holds with
\[
Q(x) = \frac{1}{2} \| x \|^2 \qquad \text{and} \qquad \gamma_1 = \gamma_2 = \gamma_4 = 1
\]
for all sufficiently large constants $\gamma_3$.
\end{example}

\begin{example}
Let $f(x) = x-x^3$, $G(x) = 0$, and the cost $\ell(x,u) = x^2 + x^4 + u^2$ be defined on the unbounded domain \mbox{$\Omega = \mathbb R$}. In this case, Assumption~\ref{ass::Q} is satisfied for $\mu = 0$,
\[
Q(x) = \frac{1}{2} \| x \|^2, \ \ \gamma_1 = \gamma_2 = \gamma_3 = 1, \quad \text{and} \quad \gamma_4 = \frac{1}{4} \; .
\]
The associated nominal nonlinear system, however, is neither controllable nor stabilizable and, consequently, does not admit a strong control Lyapunov function.
\end{example}

\bigskip
\noindent
The above examples suggest that Assumptions~\ref{ass::blanket} and~\ref{ass::Q} are satisfied for a large class of nonlinear control systems. This statement is additionally supported by a long list of examples for generic FPKs in~\cite{Bogachev2015}, where a very similar definition of Has'minski\v{i}-Lyapunov functions is used.

\subsection{Optimal Steady States}
This section discusses methods for analyzing the expected value of the stage cost $\ell$ at the optimal steady state of the FPK. As a first step towards this goal, the following lemma exploits a variant of Has'minski\v{i}'s theorem in order to guarantee the existence of a potentially suboptimal steady state and to derive an associated upper bound on the expected value of $\ell$.

\begin{lemma}
\label{lem::STATE}
Let Assumptions~\ref{ass::blanket} and~\ref{ass::Q} be satisfied and let $\mu \in C^1(\Omega)$ be a control law for which~\eqref{eq::RelaxedLyapunov} holds. Then there exists an equilibrium $\rho_\mathrm{s} \in H^1(\Omega) \cap L^1(\Omega)$ such that
\begin{eqnarray}
\label{eq::rhos}
\begin{array}{rcll}
0 & \ = \ & \mathcal A \rho_\mathrm{s} + \mathcal B(\mu \rho_\mathrm{s}) \ & \mathrm{on} \ \ \Omega, \\[0.16cm]
0 & \ \geq \ & \rho_\mathrm{s} \ & \mathrm{on} \ \ \Omega, \\[0.16cm]
0 & \ = \ & \left( \epsilon \nabla \rho_\mathrm{s} - F[\mu] \rho_\mathrm{s} \right)^\tr n \quad & \mathrm{on} \ \ \partial \Omega, \\[0.16cm]
1 & \ = \ & \int_{\Omega} \rho_\mathrm{s}(x) \, \mathrm{d}x \; , & \\[0.16cm]
\dfrac{\epsilon \gamma_2 n_x + \gamma_3}{\gamma_4} & \ \geq \  & \int_{\Omega} \ell(x,\mu(x)) \rho_\mathrm{s}(x) \, \mathrm{d}x  \; .
\end{array}
\end{eqnarray}
\end{lemma}

\bigskip
\noindent
\textbf{Proof.} Let $Q$ be defined as in Assumption~\ref{ass::Q} and let
\[
\Omega_{\alpha} \ \defeq \ \{ x \in \Omega \mid Q(x) < \alpha \} \qquad \text{for} \qquad \alpha \ > \ \inf_{x \in \Omega} Q(x)
\]
denote the associated open $\alpha$-sublevel sets of $Q$. Because $Q$ is strongly convex and differentiable and because $\alpha$ is assumed to be sufficiently large, the sets $\Omega_\alpha$ are open, non-empty, convex, bounded, and have a Lipschitz boundary. Thus, it makes sense to introduce the convex PDE-constrained auxiliary optimization problem
\begin{eqnarray}
\label{eq::linfaux}
\widetilde \ell_{\alpha} \ \defeq \ &\sup_{\rho}& \ \int_{\Omega_{\alpha}} \gamma_4 \ell[\mu] \rho \, \mathrm{d}x \\[0.2cm]
&\mathrm{s.t.}& \ \left\{
\begin{array}{rcll}
0 &=& \epsilon \Delta \rho - \mathrm{div}(F[\mu] \rho) \ \ & \text{on} \ \Omega_{\alpha} \\[0.16cm]
0 &\geq& \rho & \text{on} \ \Omega_{\alpha} \; \\[0.16cm]
0 & \ = \ & \left( \epsilon \nabla \rho - F[\mu] \rho \right)^\tr n_{\alpha} \quad & \text{on} \ \ \partial \Omega_{\alpha}, \\[0.16cm]
1 &=& \int_{\Omega_{\alpha}} \rho \, \mathrm{d}x 
\end{array}
\right. \notag
\end{eqnarray}
with optimization variable $\rho \in H^1(\Omega_\alpha) \cap L^1(\Omega_\alpha)$. Here, $n_{\alpha}$ denotes an outer normal of $\Omega_{\alpha}$. Note that the generalized version of Has'minski\v{i}'s theorem on the existence of unique probability solutions of the stationary FPK on convex domains~\cite[Thm.~A]{Huang2015}, see also~\cite{Bogachev2015}, implies that~\eqref{eq::linfaux} is feasible. Since we assume that $\ell$ is bounded from below, it follows that there exists a $\underline{\sigma} > -\infty$ with
$\underline{\sigma} < \widetilde \ell_\alpha$. Moreover, a Lagrangian relaxation of~\eqref{eq::linfaux} yields
\begin{align}
& \widetilde \ell_{\alpha} \ \leq \ \sup_{\rho \geq 0} \ \left\{ \overline{\sigma} + (\gamma_4 \ell[\mu] + F[\mu]^\tr \nabla Q - \overline{\sigma},\rho)_{L^2(\Omega_\alpha)} \right. \notag \\
& \left. \qquad \qquad \qquad - \epsilon (\nabla Q^\tr, \nabla \rho)_{L^2(\Omega_\alpha)} \right\} \notag \\[0.16cm]
& \ = \ \sup_{\rho \geq 0} \ \left\{ \vphantom{\int_{\partial \Omega}} \overline{\sigma} + (\gamma_4 \ell[\mu] + F[\mu]^\tr \nabla Q + \epsilon \Delta Q - \overline{\sigma},\rho)_{L^2(\Omega_\alpha)} \right. \notag \\
& \left. \qquad \qquad \quad - \epsilon \int_{\partial \Omega} \rho \nabla Q^\tr n_{\alpha} \mathrm{d}x \right\} \ \leq \ \overline{\sigma} \; , \label{eq::ellupper}
\end{align}
for sufficiently large constants $\overline \sigma > 0$. The latter inequality follows from~\eqref{eq::HessianBound},~\eqref{eq::RelaxedLyapunov}, and the fact that
\[
\nabla Q^\tr n \geq 0 \ \ \text{on} \ \ \partial \Omega \qquad \Longrightarrow \qquad \nabla Q^\tr n_{\alpha} \geq 0 \ \ \text{on} \ \ \partial \Omega_{\alpha}
\]
holds due to the above construction of $\Omega_{\alpha}$. More precisely,~\eqref{eq::HessianBound} and~\eqref{eq::RelaxedLyapunov} imply that
\[
\Delta Q \ \leq \ n_x \gamma_2 \quad \text{and} \quad \gamma_4 \ell[\mu] + F[\mu]^\tr \nabla Q \ \leq \ \gamma_3 \; .
\]
Consequently,~\eqref{eq::ellupper} holds for
\[
\overline \sigma \ = \ \gamma_3 + \epsilon \, n_x \gamma_2 \; .
\]
But this particular choice for $\overline \sigma$ does not depend on $\alpha$ and, consequently, the bounds $\underline \sigma \leq \ell_{\alpha} \leq \overline \sigma$ hold uniformly. The statement of the theorem follows then by taking the limit $\alpha \to \infty$ in~\eqref{eq::linfaux} and dividing the upper bound on $\widetilde \ell_\infty$ by $\gamma_4$.
\qed

\bigskip
\noindent
The stage cost of~\eqref{eq::EOC} can be represented by the convex perspective function
\begin{align}
\label{eq::ell}
\mathfrak L(\rho,\nu) \ \defeq \ \int_{\Omega} \ell \left(\cdot, \frac{\nu}{\rho} \right) \rho \, \mathrm{d}x .
\end{align}
It is defined for any strictly positive and measurable function  $\rho > 0$ and any measurable function $\nu$ for which the above integral exists. Similarly, control constraints are represented by the convex cone
\[
\mathcal K \ \defeq \ \left\{
\ (\rho,\nu) \in H^1(\Omega) \times L^2(\Omega)^{n_u} \ \middle|
\begin{array}{ll}
\nu \in \rho \mathbb U \quad & \text{on} \ \Omega \\[0.16cm]
\rho \geq 0 & \text{on} \ \Omega
\end{array}
\right\}.
\]
Thus, optimal steady states can be found by solving the convex optimization problem 
\begin{eqnarray}
\label{eq::linf}
\ell_\infty \defeq &\min_{(\rho,\nu) \in \mathcal K}& \ \mathfrak L(\rho,\nu) \\[0.2cm]
&\mathrm{s.t.}& \ \left\{
\begin{array}{rcll}
0 &=& \mathcal A \rho + \mathcal B \nu \ \ & \text{on} \ \Omega \\[0.16cm]
0 & \ = \ & \left( \epsilon \nabla \rho - f \rho - G \nu \right)^\tr n \quad & \text{on} \ \ \partial \Omega, \\[0.16cm]
1 &=& \int_{\Omega} \rho \, \mathrm{d}x ,
\end{array}
\right. \notag
\end{eqnarray}
with $(\rho,\nu) \in \mathcal K$ denoting optimization variables.

\begin{corollary}
\label{cor::minimizer}
Let Assumptions~\ref{ass::blanket} and~\ref{ass::Q} be satisfied. Then~\eqref{eq::linf} has a unique minimizer
$(\rho_{\infty},\nu_\infty) \in \mathcal K $.
\end{corollary}

\textbf{Proof.} The construction of $\rho_\mathrm{s}$ and $\mu$ in Lemma~\ref{lem::STATE} satisfies
\[
(\rho_\mathrm{s},\mu \rho_\mathrm{s}) \ \in \ \mathrm{int}(\mathcal K) \; .
\]
As such, the closed and convex cone $\mathcal K$ is compatible with Slater's constraint qualification, recalling that $\rho_\mathrm{s} > 0$ follows from the minimum principle for elliptic operators~\cite{Bogachev2015}. Consequently, feasibility and boundedness of~\eqref{eq::linf} follows immediately from Lemma~\ref{lem::STATE}, as $(\rho,\nu) = (\rho_\mathrm{s},\mu \rho_\mathrm{s})$ is a feasible point with bounded objective value. In summary,~\eqref{eq::linf} has the following properties:
\begin{enumerate}

\item Slater's (weak) constraint qualification holds,

\item the objective function $\mathcal L$ is a convex perspective function that is---due to Assumption~\ref{ass::blanket}---strongly convex in the control $\nu$, and

\item the operator $\mathcal A$ is uniformly elliptic.

\end{enumerate}

\bigskip
\noindent
Consequently, the statement of this corollary follows from a well-known standard result on the existence and uniqueness of minimizers of elliptic PDE-constrained convex optimization problems in Hilbert-Sobolev spaces; see, for example,~\cite[Thm~1.45]{Hinze2009}.\qed

\subsection{Dissipation of Energy}
Let $(\rho_\infty,\nu_\infty)$ denote the optimal steady state, as defined in Corollary~\ref{cor::minimizer}, and let $\mu_\infty = \nu_\infty/\rho_\infty$ denote an associated optimal control law. It can be used to study the controlled parabolic FPK
\begin{eqnarray}
\label{eq::FPinf}
\begin{array}{rcll}
\dot \rho & \ = \ & \mathcal A \rho + \mathcal B( \mu_\infty \rho ) & \text{on} \ \ \Omega_{\mathbb T} \\[0.16cm]
\rho_0 &=& \rho(0) & \text{on} \ \ \Omega \\[0.16cm]
0 &=& \left. \left( \epsilon \nabla \rho - F[\mu_\infty]^\tr \rho \right)^\tr n \right|_{\partial \Omega} \quad & \text{on} \ \ \Sigma_{\mathbb T}
\end{array}
\end{eqnarray}
on a possibly infinite time horizon $\mathbb T$. Here, $\rho_0 \in \mathcal P_0(\rho_\infty^{-1})$ denotes an initial probability density that is assumed to be an element of the set
\[
\mathcal P_0(\rho_\infty^{-1}) \ \defeq \ \left\{ \ \rho_0 \in L^2(\rho_\infty^{-1}) \cap L^1(\Omega) \ \middle|
\begin{array}{l}
\int_{\Omega} \rho_0 \, \mathrm{d}x = 1, \\[0.16cm]
\rho_0 \geq 0 \ \ \text{on} \ \ \Omega 
\end{array}
\right\}.
\]
In this setting, the key idea for analyzing~\eqref{eq::FPinf} is to introduce the energy functional
\[
E \ \defeq \ \| \rho - \rho_\infty \|_{L^2(\rho_\infty^{-1})}^2 \ = \ \int_{\Omega} \frac{(\rho-\rho_\infty)^2}{\rho_\infty} \, \mathrm{d}x \ . \
\]
Its weak time derivative is given by
\begin{align}
\label{eq::dissipation}
\frac{\mathrm{d}E}{\mathrm{d}t} \ = \ -2 \epsilon \int_{\Omega} \left\| \nabla \left( \frac{\rho}{\rho_\infty} \right) \right\|_2^2 \rho_\infty \, \mathrm{d}x \ \leq \ \ 0 \; .
\end{align}
Note that this equation is well-known in the functional analysis and FPK literature~\cite{Arnold2001,Bogachev2015}, where $E$ is often interpreted as energy (or entropy) that is dissipated during the evolution of the FPK. Alternatively, by translating this to control engineering language, $E$ could also be called a Lyapunov function for~\eqref{eq::FPinf}.

\begin{lemma}
\label{lem::gamma}
Let Assumptions~\ref{ass::blanket} and~\ref{ass::Q} be satisfied. Then,~\eqref{eq::FPinf} has for every $\rho_0 \in \mathcal P_0(\rho_\infty^{-1})$ a unique solution
\[
\rho \in \mathcal P(\rho_\infty^{-1}) \qquad \text{with} \qquad \left\{
\begin{array}{ll}
\rho > 0 \ & \mathrm{on} \ \ \Omega_{\mathbb T} \\[0.16cm]
\int_{\Omega} \rho \, \mathrm{d}x = 1 \ & \mathrm{on} \ \ \mathbb T.
\end{array}
\right.
\]
Moreover, $\rho$ satisfies~\eqref{eq::dissipation}.
This results holds on finite and infinite time horizons $\mathbb T$.
\end{lemma}

\textbf{Proof.} FPKs have been analyzed by so many authors that is impossible to list them all. As such,~\cite{Arnold2001} and \cite{Bogachev2015} should be regarded as selected examples for two out of a long list of articles in which proofs of variants of the current lemma can be found. Nevertheless, because this is important for the developments in this paper, it appears appropriate to recall the following three main arguments why the above statement holds.

\begin{enumerate}
\addtolength{\itemsep}{2pt}

\item First, if $\Omega_{\mathbb T}$ is bounded, a unique weak solution $\rho \in \mathcal P(\rho_\infty^{-1})$ of~\eqref{eq::FPinf} exists, since it is a uniformly parabolic linear PDE. Moreover, it follows from Harnack's inequality that
\begin{align}
\label{eq::phi}
\phi \ \defeq \ \frac{\rho}{\rho_\infty} \ \in \ H^1(\rho_\infty) \qquad \text{on} \ \ \mathbb T .
\end{align}
The details of this argument can be found in~\cite[Chapter~7]{Bogachev2015}, where the square-integrability of logarithmic gradients of FPKs is discussed in all detail.

\item Next,~\eqref{eq::dissipation} can be verified by direct computation,
\begin{eqnarray}
\frac{\mathrm{d}E}{\mathrm{d}t} & \ = \ & 2 \left\langle \dot \rho ,  \phi - 1 \right\rangle \ = \ 2 \left\langle \dot \rho ,  \phi - 1 \right\rangle - \left\langle \dot \rho_\infty ,  \phi^2 \right\rangle \notag \\[0.16cm]
&\overset{\eqref{eq::FPinf}}{=}& 2 \int_{\Omega} \left( \rho F[\mu_\infty]^\tr \nabla \phi - \epsilon \nabla \rho^\tr \nabla \phi \right) \, \mathrm{d}x \notag \\[0.16cm]
& & -2 \int_{\Omega} \left[ \rho_\infty F[\mu_\infty]^\tr \nabla \phi - \epsilon \nabla \rho_\infty^\tr \nabla \phi \right] \phi \, \mathrm{d}x
\notag \\[0.16cm]
&\overset{\eqref{eq::phi}}{=}& - 2 \epsilon \int_{\Omega} \left\| \nabla \phi \right\|_2^2 \rho_\infty \, \mathrm{d}x \ \leq \ 0,
\end{eqnarray}
where the equation in the second line follows by using that the steady state $\rho_\infty$ also satisfies~\eqref{eq::FPinf} but with $\rho_0 = \rho_\infty$ and $\dot \rho_\infty = 0$. Thus, Gronwall's lemma yields the a-priori energy estimate
\begin{align}
\label{eq::energy}
\| \rho(t) - \rho_\infty \|_{L^2(\rho_\infty^{-1})} \ \leq \ \| \rho_0 - \rho_\infty \|_{L^2(\rho_\infty^{-1})} \quad \text{on} \ \ \mathbb T.
\end{align}

\item For the case that $\Omega$ and $\mathbb T$ are unbounded, one can proceed as in the proof of Lemma~\ref{lem::STATE}. This means that the above two steps are repeated in order to construct a sequence of solutions $\rho_k \in \mathcal P(\rho_\infty^{-1})$ of~\eqref{eq::FPinf} on the bounded domains $\Omega_{\alpha_k} \times (0,T_k)$ with increasing sequences $\alpha_k \to \infty$ and $T_k \to \infty$. Because the energy estimate~\eqref{eq::energy} yields a uniform upper bound, the sequence $(\rho_k)_{k \in \mathbb N}$ has a weakly convergent subsequence in the Hilbert space $\mathcal P(\rho_\infty^{-1})$, whose limit must be a weak solution of~\eqref{eq::FPinf} and, by continuity,~\eqref{eq::dissipation} holds, too.
\end{enumerate}

\bigskip
\noindent
In addition to the above three arguments, it is also recalled that FPKs model the evolution of probability densities and, as such, the conservation law
\[
\int_{\Omega} \rho \, \mathrm{d}x \ = \ \int_{\Omega} \rho_0 \, \mathrm{d}x \ = \ 1
\]
holds on $\mathbb T$, referring to~\cite{Bogachev2015} for a formal discussion of this property. This completes the proof of the lemma. \qed

\begin{remark}
Note that more general strategies for constructing Lyapunov functions for FPKs can be found in~\cite{Arnold2001,Bolley2010}, where so-called $\Phi$-entropies are proposed as Lyapunov functions. These functions are useful for analyzing the properties of FPKs in general $L^p$-spaces, typically for $1 \leq p < \infty$. Here, the case $p=1$ is particularly relevant, as the $L^1$-norm is  a natural choice for measuring the distance between probability densities. The methods in this paper are nevertheless presented in Hilbert-Sobolev spaces, as many of the constructions below rely on the concept of duality in convex PDE-constrained optimization. Extensions to other reflexive and non-reflexive $L^p$-spaces are beyond the scope of this article.
\end{remark}

\subsection{Ergodic Optimal Control}
In order to analyze the ergodic optimal control problem~\eqref{eq::EOC}, an additional assumption is introduced.
\begin{assumption}
\label{ass::VAR}
We have $\ell[\mu_\infty] \in L^2(\rho_\infty)$.
\end{assumption}
This assumption is equivalent to requiring that $\ell[\mu_\infty]$ has a variance at the steady state $\rho_\infty$; that is,
\[
\int_{\Omega} \ell[\mu_\infty]^2 \rho_\infty \, \mathrm{d}x - \ell_\infty^2 \ < \ \infty.
\]
If Assumption~\ref{ass::blanket} and~\ref{ass::Q} hold and if $\Omega$ is bounded, this is clearly the case, and Assumption~\ref{ass::VAR} automatically holds. Otherwise, if $\Omega$ is unbounded, such a variance exists under very mild growth conditions on $\ell$, as discussed in all detail in~\cite[Chapter~3]{Bogachev2015}, where upper bound estimates for $\rho_\infty$ can be found. Next, let
\begin{eqnarray}
\label{eq::ocp}
J(T,\rho_0) \defeq &\min_{(\rho,\nu) \in \mathcal K_{\mathbb T}}& \ \int_{0}^T \mathfrak L(\rho,\nu) \, \mathrm{d}t \\[0.2cm]
&\mathrm{s.t.}& \hspace{-0.2cm} \left\{
\begin{array}{rcll}
\dot \rho &=& \mathcal A \rho + \mathcal B \nu & \text{on} \ \Omega_{\mathbb T} \\[0.16cm]
0 &=& \left( \epsilon \nabla \rho - f \rho - G \nu \right)^\tr n \ & \text{on} \ \Sigma_{\mathbb T} \\[0.16cm]
\rho_0 &=& \rho(0), & \text{on} \ \Omega
\end{array}
\right. \notag
\end{eqnarray}
denote the value function of the stochastic control problem on the finite time horizon $\mathbb T = (0,T)$ with parameters $T > 0$ and $\rho_0 \in \mathcal P_0(\rho_\infty^{-1})$. Here, the optimization variables $\rho$ and $\nu$ are elements of the convex cone
\[
\mathcal K_{\mathbb T} \ \defeq \ \left\{
\ (\rho,\nu) \ \middle|
\begin{array}{l}
\rho \in \mathcal P(\rho_\infty^{-1}) \\[0.16cm]
\nu \in \mathcal U(\rho_\infty^{-1}) \\[0.16cm]
(\rho(t),\nu(t)) \in \mathcal K \quad \text{f.a.} \ \ t \in \mathbb T
\end{array}
\right\} \; .
\]
By construction, Problem~\eqref{eq::EOC} is equivalent to computing the limit
\begin{align}
\label{eq::EOC2}
J^\star \ = \ \lim_{T \to \infty} \ \frac{J(T,\rho_0)}{T} \; .
\end{align}
The following lemma proves that $J^\star = \ell_\infty$, which, in turn, implies that $J^\star$ does not depend on $\rho_0$.

\begin{theorem}
\label{thm::ergodic}
Let Assumptions~\ref{ass::blanket},~\ref{ass::Q}, and~\ref{ass::VAR} be satisfied and let $\ell_\infty$ be defined as in~\eqref{eq::linf}. Then~\eqref{eq::ocp} has for all $T > 0$ and all $\rho_0 \in \mathcal P(\rho_\infty^{-1})$ a unique minimizer $(\rho^\star,\nu^\star) \ \in \ \mathcal K_{\mathbb T}$. Moreover, the ergodic limit~\eqref{eq::EOC2} exists. It is given by
\[
\forall \rho_0 \in \mathcal P_0(\rho_\infty^{-1}), \qquad \lim_{T \to \infty} \frac{J(T,\rho_0)}{T} \ = \ \ell_\infty \; . 
\]
\end{theorem}

\textbf{Proof.} This proof is divided into two parts, where the first part verifies the claim that~\eqref{eq::ocp} has a unique minimizer while the second part is about the ergodic limit.

\bigskip
\noindent
\textit{Part I.} First note that an upper bound on the stage cost of~\eqref{eq::ocp} can be found by applying the triangle- and the Cauchy-Schwarz inequality,
\begin{eqnarray}
\left| \mathfrak L(\rho,\nu) \right| & \ \leq \ & \ell_\infty + \left| \mathfrak L(\rho,\nu) - \ell_\infty \right| \notag
\\[0.16cm]
& = & \ \ell_\infty + \left| \int_{\Omega} \ell[\mu_\infty] (\rho - \rho_\infty) \, \mathrm{d}x \right| \notag \\[0.16cm]
& \leq &  \ell_\infty + \| \rho - \rho_\infty \|_{L^2(\rho_\infty^{-1})}
\left\| \ell[ \mu_\infty ] \right\|_{L^2(\rho_\infty)}, \label{eq::CS}
\end{eqnarray}
which holds along the trajectory $(\rho,\nu) = (\rho, \mu_\infty \rho)$ of~\eqref{eq::FPinf}. Due to Assumption~\ref{ass::VAR} and the energy estimate~\eqref{eq::energy} from the proof of Lemma~\ref{lem::gamma}, this upper bound is sufficient to ensure that $J(T,\rho_0) < \infty$. Thus, the first part of this proof is analogous to the argument in the proof of Corollary~\ref{cor::minimizer}: the conic constraint $(\rho,\nu) \in \mathcal K_{\mathbb T}$ is strictly feasible, which implies that Slater's weak constraint qualification holds, and $\mathcal L$ is strictly convex in $\nu$. Consequently, the existence of a unique minimizer of~\eqref{eq::ocp} for $T < \infty$ follows from the standard existence theorem for PDE-constrained optimization problems~\cite[Thm.~1.45]{Hinze2009}.

\bigskip
\noindent
\textit{Part II.} Lemma~\ref{lem::gamma} implies first that the solution $\rho$ of~\eqref{eq::FPinf} is such that that the function
\[
\phi(t) = \frac{\rho(t)}{\rho_\infty} \ \ \overset{t \to \infty}{\longrightarrow} \ \ C
\]
converges to a constant $C$. But we also have
\[
\| \rho(t) \|_{L^1(\Omega)} \ = \ \| \rho_\infty(t) \|_{L^1(\Omega)} \ = \ 1,
\]
since $\rho$ and $\rho_\infty$ are both probability density functions. Consequently, this constant is given by $C = 1$ and $\rho(t)$ converges to $\rho_\infty$ in $L^2(\rho_\infty^{-1})$ for $t \to \infty$. By substituting $\nu = \mu_\infty \rho$ in~\eqref{eq::ocp}, it follows that there exists for every $\delta > 0$ a $T_{\delta} < \infty$, with $T_\delta \to \infty$ for $\delta \to 0$, such that
\begin{align}
\label{eq::upperJ}
\forall T \geq T_\delta, \qquad \frac{J(T,\rho_0)}{T} \ \leq \ \ell_\infty + \delta \; .
\end{align}
Next, let $(\rho_\delta,\nu_\delta)$ denote the minimizer of~\eqref{eq::ocp} on the time horizon $\mathbb T = (0,T_\delta)$ and let
\[
\overline \rho_\delta \ \defeq \ \frac{1}{T_\delta} \int_0^{T_\delta} \rho_\delta \, \mathrm{d}t  \qquad
\text{and} \qquad  \overline \nu_\delta \ \defeq \ \frac{1}{T_\delta} \int_0^{T_\delta} \nu_\delta \, \mathrm{d}t
\]
denote their time averages. Since $\mathcal K$ is a convex cone, it follows that $(\overline \rho_\delta, \overline \nu_\delta) \in \mathcal K$ and
\begin{align}
& \mathcal A \overline \rho_\delta + \mathcal B \overline \nu_\delta \ = \ \frac{1}{T_\delta} \int_0^{T_\delta} \dot \rho_\delta \, \mathrm{d}t \ = \ \frac{\rho(T_\delta)-\rho_0}{T_\delta}, \notag
\end{align}
which implies that
\begin{align}
\label{eq::lowerJ1}
\lim_{\delta \to 0} \, \left( \mathcal A \overline \rho_\delta + \mathcal B \overline \nu_\delta \right) \ = \ 0,
\end{align}
almost everywhere on $\Omega$, since $\ell$ is bounded from below and radially unbounded. In other words, $(\overline \rho_\delta, \overline \nu_\delta)$ converges to a feasible steady-state in $\mathcal K$. Moreover, Jensen's inequality yields
\begin{align}
\label{eq::lowerJ2}
\frac{J(T_\delta,\rho)}{T_\delta} \ = \  \frac{\int_0^{T_\delta} \mathfrak L(\rho_\delta,v_\delta) \, \mathrm{d}t}{T_\delta}  \ \geq \ \mathfrak L \left( \overline \rho_\delta, \overline \nu_\delta  \right) \, .
\end{align}
But then, on the one hand,~\eqref{eq::upperJ} implies
\[
\limsup_{\delta \to 0} \ \frac{J(T_\delta,\rho_0)}{T_\delta} \ \leq \ \ell_{\infty}
\]
and, on the other hand,~\eqref{eq::lowerJ1} and~\eqref{eq::lowerJ2} and the fact that $\mathcal A$ is uniformly elliptic imply, by continuity,
\[
\liminf_{\delta \to 0} \ \frac{J(T_\delta,\rho_0)}{T_\delta} \ \geq \ \ell_{\infty} \; .
\]
Consequently, the ergodic limit for $\delta \to 0$, or, equivalently, $T_\delta \to \infty$, exists, and $J^\star = \ell_\infty$. This completes the proof of the theorem.\qed

\section{Hamilton-Jacobi-Bellman Equations}
\label{sec::HJB}
Many classical methods for analyzing the existence of solutions to HJBs~\cite{Bardi1997,Fleming1993} start by establishing (or sometimes even imposing) certain properties---for example, H\"older continuity or differentiability---of the value function of an optimal control problem and then show that this function actually satisfies its associated HJB, be it in a classical-, weak-, or viscosity sense. An impressive exception is the work by Krylov~\cite{Krylov1987,Krylov2008}, which is using and developing methods for general nonlinear parabolic PDEs that can be applied to analyze the solutions of controlled diffusions and nonlinear HJBs. Nevertheless, it appears that the convex duality of HJBs and FPKs, which, in principle, is known since a long time~\cite{Bismut1973,Fleming1989}, has never been fully exploited for the purpose of analyzing HJBs. Therefore, rather than relying on general nonlinear analysis methods, this section proposes to analyze HJBs by using modern PDE-constrained convex optimization methods~\cite{Hinze2009}, which can exploit the above mentioned duality of HJBs and FPKs.

\subsection{Hamiltonians}
Let $H: \overline \Omega \times \mathbb R^{n_x} \to \mathbb R$ and $u^\star: \overline \Omega \times \mathbb R^{n_x} \to \mathbb U$ denote the Hamiltonian and its associated minimizer,
\begin{eqnarray}
H(x,\lambda) &\defeq& \min_{u \in \mathbb U} \left\{ \ell(x,u) + \lambda^\tr [f(x) + G(x) u] \right\} \notag \\[0.16cm]
\text{and} \quad u^\star(x,\lambda) & \ \defeq \ & \underset{u \in \mathbb U}{\text{argmin}} \left\{ \ell(x,u) + \lambda^\tr [f(x) + G(x) u] \right\} \notag
\end{eqnarray}
for all ($x,\lambda) \in \overline \Omega \times \mathbb R^{n_x}$. Assumption~\ref{ass::blanket} implies that $u^\star$ is Lipschitz continuous on $\overline \Omega \times \mathbb R^{n_x}$. Thus, $H$ is Lipschitz continuous and, consequently,
\begin{eqnarray}
\label{eq::HVL2}
\forall V \in \mathcal P(\rho_\infty), \qquad H(\cdot,\nabla V) \in L^2(0,T;L^2(\rho_\infty)) \; ,
\end{eqnarray}
recalling that $\rho_\infty$ denotes the optimal steady state, as defined in Corollary~\ref{cor::minimizer}.

\subsection{Parabolic HJBs}
On finite time horizons $\mathbb T = (0,T)$ the parabolic HJB with Neumann boundary condition is given by
\begin{eqnarray}
\label{eq::HJB}
\begin{array}{rcll}
-\dot V & \ = \ & H(x,\nabla V) + \epsilon \, \Delta V \ \ & \text{on} \ \ \Omega_{\mathbb T} \\[0.16cm]
0 & \ = \ & V(T)  & \text{on} \ \ \Omega \\[0.16cm]
0 &=& \nabla V^\tr n & \text{on} \ \ \Sigma_{\mathbb T}.
\end{array}
\end{eqnarray}
The definition of weak solutions of~\eqref{eq::HJB} relies on~\eqref{eq::HVL2}.
\begin{definition}
\label{def::HJB}
A function $V \in \mathcal P(\rho_\infty)$ with $V(T,\cdot) = 0$ is called a weak solution of~\eqref{eq::HJB} if
\begin{eqnarray}
-\int_{0}^T \langle \dot V, \rho \rangle \, \mathrm{d}t \ = \int_{0}^T \int_{\Omega} \left[ \rho  H(x,\nabla V) - \epsilon \nabla \rho^\tr \nabla V \right] \mathrm{d}x\,  \mathrm{d}t \notag
\end{eqnarray}
for all test functions $\rho \in \mathcal P(\rho_\infty^{-1})$.
\end{definition}
Different from classical methods for analyzing HJBs, the following lemma uses Riesz' representation theorem to establish the existence and uniqueness of weak solutions in the Hilbert space $\mathcal P(\rho_\infty)$.

\begin{lemma}
\label{lem::HJB}
Let Assumptions~\ref{ass::blanket},~\ref{ass::Q}, and~\ref{ass::VAR} hold and let $\rho_\infty$ be defined as in Corollary~\ref{cor::minimizer}. The statements below hold for all $T < \infty$.
\begin{enumerate}
\addtolength{\itemsep}{2pt}

\item HJB~\eqref{eq::HJB} has a unique weak solution $V \in \mathcal P(\rho_\infty)$.

\item The value function $J(T,\cdot)$ of~\eqref{eq::ocp} is a bounded linear operator, whose Riesz representation is given by
\[
\forall \rho_0 \in \mathcal P_0(\rho_\infty^{-1}), \ \ \ J(T,\rho_0) \ = \ \int_\Omega V(0,x) \rho_0(x) \, \mathrm{d}x .
\]
\end{enumerate}
\end{lemma}

\textbf{Proof.} Let us temporarily assume that $\Omega$ is bounded. In this case, the auxiliary optimization problem
\begin{eqnarray}
\label{eq::JP}
\mathfrak J(\mu) \ \defeq \ & \underset{(\rho,\nu) \in \mathcal K_{\mathbb T}}{\min}& \ \int_{0}^T \mathfrak L(\rho,\nu) \, \mathrm{d}t \\[0.2cm]
&\mathrm{s.t.}& \left\{
\begin{array}{rcll}
\dot \rho &=& \mathcal A \rho + \mathcal B \nu & \text{on} \ \Omega_{\mathbb T} \\[0.16cm]
\nu & = & \mu \rho & \text{on} \ \Omega_{\mathbb T} \\[0.16cm]
\rho_0 &=& \rho(0), & \text{on} \ \Omega \\[0.16cm]
0 &=& \left( \epsilon \nabla \rho - f \rho - G \nu \right)^\tr n \ \ & \text{on} \ \Sigma_{\mathbb T}
\end{array}
\right. \notag
\end{eqnarray}
satisfies Slater's weak constraint qualification for uniformly parabolic PDE-constrained convex optimization problems on a bounded domain, recalling that $\mathfrak L$ is strongly convex in $\nu$. Consequently, it follows from a standard existence result for PDE-constrained optimization problems~\cite[Chapter~1.5]{Hinze2009} that~\eqref{eq::JP} has for all parameters $\mu \in \mathfrak U$ and given $\rho_0 \in \mathcal P_0(\rho_\infty^{-1})$ a unique minimizer $(\rho_\mu,\nu_\mu) \in \mathcal K_{\mathbb T}$. Next, let $V \in \mathcal P(\rho_\infty)$ denote the unique weak solution of the uniformly parabolic linear PDE
\begin{eqnarray}
\label{eq::HJB1}
\begin{array}{rcll}
-\dot V &=& \ell[\mu^\star] + F[\mu^\star]^\tr \nabla V  + \epsilon \Delta V \ \ & \text{on} \ \Omega_{\mathbb T} \\[0.16cm]
0 &=& V(T) & \text{on} \ \Omega \\[0.16cm]
0 &=& \nabla V^\tr n  & \text{on} \ \Sigma_{\mathbb T}
\end{array}
\end{eqnarray}
at $\mu^\star = \nu^\star/\rho^\star$, where $(\rho^\star,\nu^\star) \in \mathcal K_{\mathbb T}$ denotes---as in Theorem~\ref{thm::ergodic}---the minimizer of~\eqref{eq::ocp}. This notation is such that the relations $\rho^\star = \rho_{\mu^\star}$ and $\nu^\star = \nu_{\mu^\star}$ hold. Let
\begin{eqnarray}
\chi(\rho,\mu) & \ \defeq \ & \int_{0}^T \{ \langle \dot V + \ell[\mu] + F[\mu]^\tr \nabla V, \rho \rangle \notag \\[0.16cm]
& & \quad \ - \epsilon \nabla \rho^\tr \nabla V \} \, \mathrm{d}t + \int_\Omega V(0,x) \rho_0(x) \, \mathrm{d}x \notag
\end{eqnarray}
denote an associated parametric Lagrangian at $V$. Because $V$ is a weak solution of~\eqref{eq::HJB1}, the function $\chi(\cdot,\mu^\star)$ is constant on $\mathcal P(\rho_\infty^{-1})$ and satisfies
\begin{eqnarray}
\label{eq::chiconst}
\forall \rho \in \mathcal P(\rho_\infty^{-1}), \quad \chi(\rho,\mu^\star) \ = \ \int_\Omega V(0,x) \rho_0(x) \, \mathrm{d}x \; .
\end{eqnarray}
Moreover, the Lagrangian $\chi$ is closely related to the Hamiltonian $H$. In order to elaborate on this relation, the auxiliary control law
\[
\hat \mu \ \defeq \ u^\star(\cdot,\nabla V) \; ,
\]
is introduced, recalling that $u^\star$ is Lipschitz continuous. Thus, it satisfies $\hat \mu \in \mathfrak U$ and
\begin{eqnarray}
\inf_{\mu \in \mathfrak U} \ \chi(\rho,\mu) & \ = \ & \chi(\rho,\hat \mu) \label{eq::chi} \\[0.16cm]
& = & \int_{0}^T \int_{\Omega} \left\{ \rho H(x,\nabla V) - \epsilon \nabla \rho^\tr \nabla V \right\}  \mathrm{d}x \, \mathrm{d}t \notag \\[0.16cm]
& & + \int_{0}^T \langle \dot V, \rho \rangle \, \mathrm{d}t + \int_\Omega V(0,x) \rho_0(x) \, \mathrm{d}x \notag
\end{eqnarray}
for all $\rho \in \mathcal P(\rho_\infty^{-1})$ with $\rho \geq 0$. Moreover, by regarding $V$ as a test function, it follows that
\begin{eqnarray}
\mathfrak J(\mu) & \ = \ & \int_{0}^T \left\{ \mathfrak L(\rho_\mu,\mu \rho_\mu) - \langle \dot \rho_\mu, V \rangle \right. \notag \\[0.16cm]
& & \qquad \qquad \qquad \quad \left. + a(\rho_\mu,V) + b(\mu \rho_\mu,V) \right\} \, \mathrm{d}t \notag \\[0.16cm]
&=& \chi(\rho_\mu,\mu) \; , \label{eq::dualJP}
\end{eqnarray}
for all $\mu \in \mathfrak U$, where the second equation in~\eqref{eq::dualJP} follows by substituting the definitions of $a,b,\chi$, and $\mathfrak L$ and a partial integration with respect to time. Next, an important observation is that
\begin{eqnarray}
\mathfrak J(\mu^\star) & \ \overset{\eqref{eq::dualJP}}{=} \ & \chi( \rho^\star, \mu^\star ) \ \overset{\eqref{eq::chiconst}}{=} \ \chi( \rho_{\hat \mu}, \mu^\star ) \ \geq \ \inf_{\mu \in \mathfrak U} \chi(\rho_{\hat \mu}, \mu ) \notag \\[0.16cm]
& \ \overset{\eqref{eq::chi}}{=} \ & \chi(\rho_{\hat \mu},\hat \mu) \ \overset{\eqref{eq::dualJP}}{=} \ \mathfrak J(\hat \mu) \; . \label{eq::Jineq}
\end{eqnarray}
Since $\mu^\star = \nu^\star/\rho^\star$ is---due to Theorem~\ref{thm::ergodic}---the unique minimizer of the optimization problem
\[
\mathfrak J(\mu^\star) \ = \ \min_{\mu \in \mathfrak U} \ \mathfrak J(\mu) \ \leq \ \mathfrak J(\hat \mu) \quad \overset{\eqref{eq::Jineq}}{\Longrightarrow} \quad \mathfrak J(\mu^\star) \ = \ J(\hat \mu),
\]
it follows from~\eqref{eq::Jineq} that $\hat \mu = \mu^\star$. Moreover, by substituting the equation $\mu^\star = \hat \mu$ in~\eqref{eq::HJB1}, we find that $V$ satisfies~\eqref{eq::HJB}; and, finally,~\eqref{eq::chiconst} and~\eqref{eq::dualJP} yield
\[
\quad J(T,\rho_0) \ = \ \mathfrak J(\mu^\star) \ = \ \int_\Omega V(0,x) \rho_0(x) \, \mathrm{d}x
\]
for all $\rho_0 \in \mathcal P(\rho_\infty^{-1})$. This implies that $J(T,\cdot)$ is a bounded linear operator. Next, if HJB~\eqref{eq::HJB} had a second solution $\widetilde V \neq V$, we could repeat the above construction by using the control law
\[
\widetilde \mu \ = \ u^\star(\cdot,\nabla \widetilde V)
\]
in order to find yet another representation of the operator $J(T,\cdot)$, namely,
\[
J(T,\rho_0) \ = \ \int_\Omega \widetilde V(0,x) \rho_0(x) \, \mathrm{d}x
\]
for any $\rho_0 \in \mathcal P(\rho_\infty^{-1})$. But this is impossible, as the Riesz representation of the bounded linear operator $J(T,\cdot)$ must be unique. Thus, HJB~\eqref{eq::HJB} has a unique solution and we have proven both statements of the theorem under the additional assumption that $\Omega$ is bounded.

\bigskip
\noindent
Last but not least, if $\Omega$ is unbounded, the same strategy as in the proof of Lemma~\ref{lem::gamma} can be applied. This means that we repeat the above argument by replacing the domain $\Omega$ in~\eqref{eq::ocp} with an increasing sequence of bounded domains $\Omega_{\alpha_k}$ in order to construct an associated sequence of bounded linear operators, $J_k(T,\cdot)$, that converges to $J(T,\cdot)$ for $k \to \infty$. Since Theorem~\ref{thm::ergodic} has already established the fact that $J(T,\cdot)$ is bounded for finite $T$ but potentially unbounded domains $\Omega$, it follows that $J(T,\cdot)$ is a bounded linear operator. Note that this is equivalent to the convergence of the corresponding sequence of Riesz' representations $V_k \in \mathcal P(\rho_\infty)$ of the operators
\[
J_k(T,\rho_0) \ = \ \int_{\Omega_{\alpha_k}} V_k(0,x) \rho_0(x) \, \mathrm{d}x,
\]
since in Hilbert spaces one can rely on the natural isometric isomorphism between linear operators and their Riesz' representations. Because $H$ is Lipschitz continuous, the limit of the sequence $V_k$ for $k \to \infty$ is the unique weak solution of~\eqref{eq::HJB}. This completes the proof of the lemma.
\qed

\begin{remark} As already mentioned above, parabolic HJBs are among the most well-studied nonlinear PDEs in the literature~\cite{Bardi1997,Fleming1993,Krylov2008}. Obviously, Lemma~\ref{lem::HJB}  adds yet another set of assumptions under which such HJBs admit a unique solution. The relevance of this result, however, is that---different from all other methods for analyzing the existence of solutions of HJBs---the proof of Lemma~\ref{lem::HJB} starts from the convex optimization problem~\eqref{eq::ocp} and then derives the HJB~\eqref{eq::HJB} by constructing an associated dual problem. As discussed in the following sections, this convex optimization based derivation of HJBs turns out to have many advantages, as it enables us to analyze HJBs by leveraging on modern functional analysis methods for linear FPKs. 
\end{remark}

\section{Exponential Stabilizability}
\label{sec::BEC}
\textit{Carr\'e Du Champ} (CDC) operators have originally been introduced in the French literature on functional analysis~\cite{Bakry1985}. They can be used to analyze general diffusions. The following section introduces a variant of such CDC operators in order to establish global exponential stabilizability conditions for controlled FPKs.

\subsection{Weighted Carr\'e Du Champ Operators}
The weighted CDC operator $\Gamma$ is a bilinear form,
\begin{align}
\label{eq::Gamma}
\forall \varphi,\phi \in C^\infty(\Omega), \qquad \Gamma(\varphi,\phi) \ \defeq \ \nabla \varphi^\tr P \nabla \phi \; .
\end{align}
Here, $P \in C^2(\Omega)^{n_x \times n_x}$ denotes a weight that satisfies $\underline \lambda I \preceq P \preceq \overline \lambda  I$ for given positive lower- and upper bounds, $0 < \underline \lambda \leq \overline \lambda < \infty$. For the special case that $P = \epsilon I$ is constant and equal to the scaled unit matrix $\epsilon I$, $\Gamma$ coincides with the traditional CDC operator of the diffusion semigroup that is generated by the linear operator
\[
\forall \phi \in C^\infty(\Omega), \qquad \mathcal L[\mu] \phi \ \defeq \ F[\mu]^\tr \nabla \phi + \epsilon \Delta \phi.
\]
It is defined for all control laws $\mu \in \mathfrak U$, although the right-hand of the CDC inequality
\begin{eqnarray}
\label{eq::CDC}
\frac{1}{2} \left( \mathcal L[\mu] \phi^2 - 2 \phi \mathcal L[\mu] \phi\right) \ = \ \epsilon \left\| \nabla \phi \right\|_2^2 \ \leq \ \frac{\epsilon}{\underline \lambda} \, \Gamma(\phi,\phi)
\end{eqnarray}
is an invariant, in the sense that it does not depend on $\mu$.

\subsection{Lyapunov Diffusion Operators}
The design of potentially non-trivial weight functions $P$ turns out to be important for the construction of practical stabilizability conditions. In order to explain why this is so, the Lyapunov diffusion operator
\begin{align}
\mathcal R[\mu] P \ \defeq \ \frac{1}{2} \left( \mathcal L[\mu] P - F'[\mu]P - PF'[\mu]^\tr \right),
\end{align}
is introduced for arbitrary $\mu \in \mathfrak U$. Here, $\mathcal L[\mu]$ is applied to $P \in C^2(\Omega)^{n_x \times n_x}$ componentwise and $F'[\mu]$ denotes the distributional Jacobian of $F[\mu]$. This means that we define
\begin{eqnarray}
(\mathcal L[\mu] P)_{i,j} & \ \defeq \ & F[\mu]^\tr \nabla P_{i,j} + \epsilon \Delta P_{i,j} \notag \\[0.16cm]
\text{and} \qquad (F'[\mu])_{i,j} & \defeq & \partial_{x_j} F_i[\mu] \notag
\end{eqnarray}
for all components $i,j \in \{ 1, \ldots, n_x \}$, where $\partial_{x_j}$ denotes the distributional derivative with respect to $x_j$. In the following, we use the convention that a formal symmetric matrix-valued distribution $Y$ on $\Omega$ is positive semi-definite, writing $Y \succeq 0$, if
\[
\int_{\Omega} \mathrm{Tr}(YZ) \, \mathrm{d}x \ \geq \ 0
\]
for all symmetric and positive semi-definite test functions $Z \in C_0^\infty(\Omega)^{n_Y \times n_Y}$ with compatible dimension. Linear matrix inequalities (LMIs) are then understood in an analogous manner.

\begin{assumption}
\label{ass::LMI}
Let $\mu_\infty = \nu_\infty/\rho_\infty$ denote the optimal ergodic control law as defined in Corollary~\ref{cor::minimizer}. We assume that there exists a $\lambda > 0$ such that the LMI
\begin{align}
\label{eq::LMI}
\left(
\begin{array}{cc}
\mathcal R[\mu_\infty] P & \nabla^\tr \otimes (\epsilon P) \\[0.16cm]
\nabla \otimes (\epsilon P) & I \otimes (\epsilon P)
\end{array}
\right) \ \succeq \
\lambda
\left(
\begin{array}{cc}
P & \ 0 \ \\
0 & 0
\end{array}
\right) \ \ \text{on} \ \ \Omega \ \
\end{align}
admits a positive definite solution $P \in C^2(\Omega)^{n_x \times n_x}$ with $\underline \lambda I \preceq P \preceq \overline \lambda P$ on $\Omega$. Here, $\otimes$ denotes the Kronecker product.
\end{assumption}

\bigskip
\noindent
Examples for linear and nonlinear control systems for which Assumption~\ref{ass::LMI} holds can be found below.

\begin{example}
Let $f(x) = Ax$ and $G(x) = B$ be a linear control system on $\Omega = \mathbb R^{n_x}$ with $(A,B)$ being asymptotically stabilizable and $\ell(x,u) = \| x \|^2 + \| u \|^2$. In this case, $\mu_\infty = K x$ is the well-known Linear-Quadratic-Gaussian (LQG) regulator and Assumption~\ref{ass::LMI} holds for a constant $P \in \mathbb R^{n_x \times n_x}$. This is because for a constant $P \succ 0$, the equations $\mathcal L[\mu] P = 0$ and $\nabla \otimes P = 0$ hold such that LMI~\eqref{eq::LMI} reduces to
\[
\mathcal R[\mu]P \ = \ -\frac{(A + B K)P + P(A+BK)^\tr}{2} \ \succeq \ \lambda P \; .
\]
This LMI has a strictly positive solution $P \succ 0$ for sufficiently small $\lambda > 0$, since we assume that $(A,B)$ is asymptotically stabilizable, or, equivalently, that the eigenvalues of $A+BK$ have strictly negative real parts.
\end{example}

\bigskip
\noindent
\begin{example}
Let us consider a stochastic system with $f(x) = \frac{x}{2} - x^3$, $G(x) = 0$, and $\epsilon = 1$ on $\Omega = \mathbb R$. It leads to an SDE with multi-modal limit distribution,
\[
\rho_\infty(x) = \frac{1}{C} \exp \left( \frac{x^2-x^4}{4} \right),
\]
where the constant $C$ can be chosen such that the integral over $\rho_\infty$ is equal to $1$. In this example,~\eqref{eq::LMI} cannot be satisfies for a constant function $P$, but it is possible to find other classes of functions which satisfy~\eqref{eq::LMI}. For instance,  Assumption~\ref{ass::LMI} can be satisfied by setting
\[
P(x) = 1 - \frac{1}{2} e^{-x^2} , \ \ \lambda = \frac{1}{20}, \ \ \underline \lambda = \frac{1}{2}, \ \ \text{and} \ \ \overline \lambda = 1 \; .
\]
\end{example}

\bigskip
\noindent
Apart from the above examples, another reason why Assumption~\ref{ass::LMI} appears appropriate for our purposes, is elaborated in the following section. Namely, it turns out that LMI~\eqref{eq::LMI} is satisfied for an even broader class of problems than the original Bakry-Emery condition for FPKs, which is---at least in the FPK literature---well-accepted as a useful and practical condition; see, for instance,~\cite{Arnold2001,Bakry1985,Bolley2010,Ledoux1992}, where many examples are collected.

\subsection{Exponential Convergence Rate}

The following lemma interprets LMI~\eqref{eq::LMI} as a generalized Bakry-Emery condition.

\begin{lemma}
Let Assumptions~\ref{ass::blanket},~\ref{ass::Q}, and~\ref{ass::LMI} hold and let \mbox{$\rho_0 \in \mathcal P_0(\rho_\infty^{-1})$} be given. Then the controlled FPK
\begin{eqnarray}
\label{eq::FokkerPlanckAUX}
\begin{array}{rcll}
\dot \rho &=& \mathcal A \rho + \mathcal B (\mu_\infty \rho) & \text{on} \ \ \Omega_{\mathbb T} \\[0.16cm]
\rho_0 &=& \rho(0) & \text{on} \ \ \Omega \\[0.16cm]
0 &=& \left. \left( \nabla \rho - F[\mu_\infty] \rho  \right)^\tr n \right|_{\partial \Omega} \quad & \text{on} \ \ \Sigma_{\mathbb T} \\
\end{array}
\end{eqnarray}
has a unique solution $\rho \in \mathcal P(\rho_\infty^{-1})$ that converges to the optimal ergodic density function $\rho_\infty$ in $L^2(\rho_\infty^{-1})$ with exponential convergence rate
$\gamma = 2 \lambda \cdot \left( \underline \lambda / \overline \lambda \right)$; that is,
\begin{eqnarray}
\label{eq::expconv}
\| \rho(t) - \rho_\infty \|_{L^2(\rho_\infty^{-1})} \ \leq \ e^{-\gamma t} \| \rho_0 - \rho_\infty \|_{L^2(\rho_\infty^{-1})}
\end{eqnarray}
for all $t > 0$.
\end{lemma}

\textbf{Proof.} In order to interpret LMI~\eqref{eq::LMI} as a generalized Bakry-Emery condition, we work out the $\Gamma_2$-operator
\begin{eqnarray}
\Gamma_2(\phi) & \ \defeq \ & \dfrac{1}{2} \left( \mathcal L[\mu_\infty] \Gamma(\phi,\phi) - 2 \Gamma(\phi, \mathcal L[\mu_\infty] \phi ) \right) \notag \\[0.26cm]
&=& \ \dfrac{1}{2} F[\mu_\infty]^\tr \nabla( \nabla \phi^\tr P \nabla \phi ) + \epsilon \Delta ( \nabla \phi^\tr P \nabla \phi ) \notag \\[0.16cm]
& & \hphantom{F[\mu_\infty]} - \nabla \phi P \nabla ( F[\mu_\infty]^\tr \nabla \phi + \epsilon \Delta \phi ) \notag \\[0.26cm]
&= \ & \dfrac{1}{2} \nabla \phi^\tr \mathcal L[\mu_\infty] P \nabla \phi + 2 \epsilon (\nabla \otimes \nabla \phi)^\tr (\nabla \otimes P) \nabla \phi \notag \\[0.16cm]
& & \hphantom{F[\mu_\infty]} + \epsilon (\nabla \otimes \nabla \phi)^\tr(I \otimes P)(\nabla \otimes \nabla \phi) \notag \\[0.16cm]
& & \hphantom{F[\mu]} - \nabla \phi P F'[\mu_\infty]^\tr \nabla \phi \notag \\[0.26cm]
&=& \ \left(
\begin{array}{c}
\nabla \phi \\
\nabla \otimes \nabla \phi
\end{array}
\right)^\tr
\left(
\begin{array}{cc}
\mathcal R[\mu_\infty]P & \epsilon \nabla^\tr \otimes P \\
\epsilon \nabla \otimes P & \epsilon I \otimes P
\end{array}
\right)
\left(
\begin{array}{c}
\nabla \phi \\
\nabla \otimes \nabla \phi
\end{array}
\right)
\end{eqnarray}
for arbitrary test functions $\phi \in C_0^\infty(\overline \Omega)$. Note that this equation is understood by interpreting all first order derivatives of $F[\mu_\infty]$ as distributional derivatives. Next,~\eqref{eq::LMI} implies that the integral version of the Bakry-Emery condition~\cite{Bakry1985,Ledoux1992},
\begin{align}
\label{eq::BEC}
\int_{\Omega} \Gamma_2(\phi) \omega \, \mathrm{d}x \ \geq \ \lambda \int_{\Omega} \Gamma(\phi,\phi) \omega \, \mathrm{d}x,
\end{align}
is satisfied for all test functions $\phi \in C_0^\infty(\overline \Omega)$ and all weight functions $\omega \in H^1(\Omega)$ with $\omega \geq 0$. This integral inequality remains correct even if we work with distributional Jacobians $F'[\mu_\infty]$, since these distributional derivatives merely enter via the integrand of the left-hand integral in~\eqref{eq::BEC}, and, consequently, can be resolved by partial integration. It is well-known~\cite[Thm.~2 \& Prop.~5]{Bolley2010} that~\eqref{eq::CDC} and~\eqref{eq::BEC} imply that the entropy inequality
\begin{align}
\label{eq::EntropyIneq}
\int_{\Omega} (\phi-1)^2 \, \rho_\infty \, \mathrm{d}x \ \leq \ \frac{\epsilon}{2 \lambda \underline \lambda} \int_\Omega \Gamma( \phi, \phi ) \, \rho_\infty \, \mathrm{d}x,
\end{align}
holds for all $\phi \in C_0^\infty(\overline \Omega)$, recalling that $(\rho_\infty \, \mathrm{d}x)$ is, by construction, an ergodic measure of the diffusion operator $\mathcal L[\mu_\infty]$. Of course, by continuity, it follows that~\eqref{eq::EntropyIneq} then also holds for all functions $\phi \in H^1(\rho_\infty)$. Consequently, it further follows from Lemma~\ref{lem::gamma} that
\begin{eqnarray}
\frac{\mathrm{d}}{\mathrm{d}t} \| \rho(t) - \rho_\infty \|_{L^2(\rho_\infty^{-1})}^2 & \ \leq \ & \frac{-2 \epsilon}{\overline \lambda} \int_{\Omega} \Gamma(\phi,\phi) \rho_\infty \, \mathrm{d}x \notag \\[0.16cm] 
&\overset{\eqref{eq::EntropyIneq}}{\leq}& - 4 \lambda (\underline \lambda / \overline \lambda) \int_{\Omega} (\phi-1)^2 \, \rho_\infty \, \mathrm{d}x \notag \\[0.16cm]
&=& - 2 \gamma \| \rho(t) - \rho_\infty \|_{L^2(\rho_\infty^{-1})}^2 , \label{eq::gron}
\end{eqnarray}
where we have set $\phi = \rho/\rho_\infty$. An application of Gronwall's lemma to~\eqref{eq::gron} yields the energy estimate~\eqref{eq::expconv}, which, in turn, implies that the unique weak solution $\rho \in \mathcal P(\rho_\infty^{-1})$ of~\eqref{eq::FokkerPlanckAUX} also remains well-defined on infinite horizons.\qed

\section{Infinite Horizon Optimal Control}
\label{sec::IHOCP}
The technical result of Lemma~\ref{lem::gamma} on the exponential convergence of the state of the optimally controlled FPK can be used to bound the value function of infinite horizon optimal control problems. Here, the main idea is to exploit the strong duality between FPKs and HJBs that has been established by Lemma~\ref{lem::HJB}.

\subsection{Problem Formulation}

This section concerns infinite horizon optimal control problems of the form
\begin{align}
\label{eq::Jinf}
J_\infty(\rho_0) \ \defeq \ \lim_{T \to \infty} \ \left( J(T,\rho_0) - T \ell_\infty\right) \; ,
\end{align}
recalling that the cost function $J$ has been defined in~\eqref{eq::ocp}.
Here, $\ell_\infty$ denotes the ergodic average cost that has been defined in~\eqref{eq::linf} and interpreted in Theorem~\ref{thm::ergodic}. Compared to the ergodic optimal control problem from Theorem~\ref{thm::ergodic}, the cost $J_\infty$ of the infinite horizon optimal control problem~\eqref{eq::Jinf} is more difficult to analyze, as it depends on the initial probability density function $\rho_0 \in \mathcal P_0(\rho_\infty)$. Conditions under which the limit on the right hand of~\eqref{eq::Jinf} exists are discussed below.

\subsection{Ergodic HJB}
Problem~\eqref{eq::Jinf} is closely related to the ergodic HJB with Neumann boundary condition, given by
\begin{eqnarray}
\label{eq::HJBinf}
\begin{array}{rcll}
\ell_\infty & \ = \ & H(x,\nabla V_\infty) + \epsilon \, \Delta V_\infty \ \ & \text{on} \ \ \Omega \\[0.16cm]
0 &=& \left(\nabla V_\infty\right)^\tr n & \text{on} \ \ \partial \Omega \; .
\end{array}
\end{eqnarray}
Weak solutions of~\eqref{eq::HJBinf} are defined as follows, recalling once more that $H$ is Lipschitz continuous.
\begin{definition}
\label{def::HJB2}
A function $V_\infty \in H^1(\rho_\infty)$ is called a weak solution of~\eqref{eq::HJBinf} if
\begin{eqnarray}
\ell_\infty \int_{\Omega} \rho \, \mathrm{d}x \ = \ \int_{\Omega} \left[ \rho  H(x,\nabla V_\infty) - \epsilon \nabla \rho^\tr \nabla V_\infty \right] \mathrm{d}x \notag
\end{eqnarray}
for all test functions $\rho \in H^1(\rho_\infty^{-1})$.
\end{definition}
The relation between this ergodic HJB and~\eqref{eq::Jinf} is explained in the following section.

\subsection{Solution of the Infinite Horizon Control Problem}
The following theorem summarizes a novel result about the existence of weak solutions of ergodic HJBs in the presence of general indefinite stage cost functions.

\begin{theorem}
\label{thm::HJB}
If Assumptions~\ref{ass::blanket},~\ref{ass::Q},~\ref{ass::VAR}, and~\ref{ass::LMI} hold, then the limit $J_\infty$ in~\eqref{eq::Jinf} exists for all $\rho_0 \in \mathcal P_0(\rho_\infty^{-1})$ and there exists a unique $V_\infty \in H^1(\rho_\infty)$ such that
 \[
\forall \rho_0 \in \mathcal P_0(\rho_\infty^{-1}), \qquad J_\infty(\rho_0) \ = \ \int_{\Omega} V_\infty(x) \rho_0(x) \, \mathrm{d}x \; .
 \]
Moreover, $V_\infty$ is a weak solution of the ergodic HJB~\eqref{eq::HJBinf}. It satisfies the energy estimate
\begin{align}
\label{eq::HJBenergy}
\int_{\Omega} (V_\infty)^2 \rho_\infty \, \mathrm{d}x \ \leq \ \frac{1}{\gamma^2} \int_{\Omega} \ell[\mu_\infty]^2 \rho_\infty \, \mathrm{d}x,
\end{align}
recalling that $\gamma = 2 \lambda \cdot \left( \underline \lambda / \overline \lambda \right)$. 
\end{theorem}

\bigskip
\noindent
\textbf{Proof.} The proof is divided into three parts:

\bigskip
\noindent
\textit{Part I.} Let $\rho$ denote the solution of~\eqref{eq::FokkerPlanckAUX} and let $\nu = \mu_\infty \rho$ denote its associated input function such that
\begin{eqnarray}
\left| \mathfrak L(\rho,\nu) - \ell_\infty \right| & \overset{\eqref{eq::CS}}{\leq} & \| \rho - \rho_\infty \|_{L^2(\rho_\infty^{-1})} \| \ell[\mu_\infty] \|_{L^2(\rho_\infty)} \notag \\[0.16cm]
&\overset{\eqref{eq::expconv}}{\leq}& \ e^{-\gamma t} \| \rho_0 - \rho_\infty \|_{L^2(\rho_\infty^{-1})}  \| \ell[\mu_\infty] \|_{L^2(\rho_\infty)} . \notag
\end{eqnarray}
The above upper bounded can be used to bound the value function of the auxiliary optimization problem
\begin{eqnarray}
\label{eq::JPinf}
J^+(T,\rho_0) \, \defeq \, & \underset{\rho,\nu}{\min}& \ \int_{0}^T \left( \mathfrak L(\rho,\nu) - \ell_\infty \right) \, \mathrm{d}t \\[0.2cm]
&\mathrm{s.t.}& \left\{
\begin{array}{rcll}
\dot \rho &=& \mathcal A \rho + \mathcal B \nu & \text{on} \ \Omega_{\mathbb T} \\[0.16cm]
\nu & = & \mu_\infty \rho & \text{on} \ \Omega_{\mathbb T} \\[0.16cm]
\rho_0 &=& \rho(0), & \text{on} \ \Omega \\[0.16cm]
0 &=& \left( \epsilon \nabla \rho - f \rho - G \nu \right)^\tr n \ & \text{on} \ \Sigma_{\mathbb T} .
\end{array}
\right. \notag
\end{eqnarray}
Namely, by integrating the above inequality, we find that
\begin{align}
\left| J^+(T,\rho_0) \right| \ \leq \ \frac{\| \rho_0 - \rho_\infty \|_{L^2(\rho_\infty^{-1})}}{\gamma} \sqrt{\int_{\Omega} \ell[\mu_\infty]^2 \rho_\infty \, \mathrm{d}x} \label{eq::energy1}
\end{align}
for all $T > 0$ and all $\rho_0 \in \mathcal P_0(\rho_\infty^{-1})$. Moreover, by using an analogous argument as in the proof of Lemma~\ref{lem::HJB}, we can use the weak solution of the linear parabolic PDE
\begin{eqnarray}
\begin{array}{rcll}
-\dot V^+ &=& \ell[\mu_\infty]-\ell_\infty + F[\mu_\infty]^\tr \nabla V^+  + \epsilon \Delta V^+ \ \ & \text{on} \ \Omega_{\mathbb T} \\[0.16cm]
0 &=& V^+(T) & \text{on} \ \Omega \\[0.16cm]
0 &=& \nabla (V^+)^\tr n  & \text{on} \ \Sigma_{\mathbb T}
\end{array} \notag
\end{eqnarray}
to represent $J^+(T,\rho)$ in the form
\[
\forall \rho_0 \in \mathcal P_0(\rho_\infty^{-1}), \qquad J^+(T,\rho_0) \ = \ \int_{\Omega} V^+(0,x)\rho_0(x) \, \mathrm{d}x \; .
\]
Next, because $J(T,\cdot)$ is a uniformly bounded linear operator, its limit for $T \to \infty$ admits a unique Riesz representation,
\[
J_\infty^+(\rho_0) \ \defeq \ \lim_{T \to \infty} J^+(T,\rho_0) \ = \ \int_{\Omega} V_\infty(x)\rho_0(x) \, \mathrm{d}x,
\]
where $V_\infty \in H^1(\rho_\infty)$ must be a weak solution of the linear elliptic PDE
\begin{eqnarray}
\label{eq::HJBaa}
\begin{array}{rcll}
\ell_\infty &=& \ell[\mu_\infty] + F[\mu_\infty]^\tr \nabla V_\infty + \epsilon \Delta V_\infty \ \ & \text{on} \ \Omega \\[0.16cm]
0 &=& \nabla (V_\infty)^\tr n  & \text{on} \ \partial \Omega
\end{array}
\end{eqnarray}
at $\mu_\infty = \nu_\infty/\rho_\infty$.

\bigskip
\noindent
\textit{Part II.} Let us use the function $J_\infty^+$ as a terminal cost of the finite horizon optimal control problem
\begin{eqnarray}
\widetilde J(T,\rho_0)  \defeq & \underset{(\rho,\nu) \in \mathcal K_{\mathbb T}}{\min}  & \int_{0}^T \hspace{-0.15cm} \left( \mathfrak L(\rho,\nu) - \ell_\infty \right) \, \mathrm{d}t + J_\infty^+(\rho(T)) \notag \\[0.2cm]
&\mathrm{s.t.}& \left\{
\begin{array}{rcll}
\dot \rho &=& \mathcal A \rho + \mathcal B \nu & \text{on} \ \Omega_{\mathbb T} \\[0.16cm]
\rho_0 &=& \rho(0), & \text{on} \ \Omega \\[0.16cm]
0 &=& \left( \epsilon \nabla \rho - f \rho - G \nu \right)^\tr n \ & \text{on} \ \Sigma_{\mathbb T} .
\end{array}
\right. \notag \\
\label{eq::Jtilde}
\end{eqnarray}
By using an analogous argument as in Lemma~\ref{lem::HJB}, it follows that $\widetilde J(T,\cdot)$ is a bounded linear operator, whose unique Riesz representation
\[
\widetilde J(T,\rho_0) \ = \ \int_{\Omega} \widetilde V(0,x) \rho_0(x) \, \mathrm{d}x
\]
is given by the weak solution of the parabolic HJB
\begin{eqnarray}
\label{eq::HJBtilde}
\begin{array}{rcll}
-\partial_t \widetilde V & \ = \ & H(x,\nabla \widetilde V) + \epsilon \, \Delta \widetilde V \ \ & \text{on} \ \ \Omega_{\mathbb T} \\[0.16cm]
V_\infty & \ = \ & \widetilde V(T)  & \text{on} \ \ \Omega \\[0.16cm]
0 &=& \nabla \widetilde V^\tr n & \text{on} \ \ \Sigma_{\mathbb T}.
\end{array}
\end{eqnarray}
Since $\mu_\infty$ is a feasible time-invariant control law for~\eqref{eq::Jtilde}, we must have
\[
\forall x \in \Omega, \quad \widetilde V(0,x) \ \leq \ V_\infty(x).
\]
Next, let us assume that we can find a $\rho_0 \in \mathcal P_0(\rho_\infty^{-1})$ and a $\widetilde T > 0$ for which $\widetilde J(\widetilde T,\rho_0) < J_\infty^+(\rho_0)$. Then the Lebesgue measure of the set
\[
\{ \ x \in \Omega \ \mid \ \widetilde V(0,x) \ < \ V_\infty(x) \ \}
\]
must be strictly positive. But this implies that there exists a $\delta > 0$ such that
\begin{align}
\label{eq::contradiction}
\forall T \geq \widetilde T, \qquad \widetilde J(T,\rho_\infty) \ \leq \ J_\infty^+(\rho_\infty) - \delta \ = \ -\delta \; ,
\end{align}
where $\delta$ does not depend on $T$. In order to show that this inequality leads to a contradiction, we denote the minimizer of~\eqref{eq::Jtilde} at $\rho_0 = \rho_\infty$ by $(\rho_T^\star,\nu_T^\star)$. Due to the above inequality, the ergodic average values
\[
\overline \rho \ \defeq \ \lim_{T \to \infty} \frac{\int_0^T \rho_T^\star \, \mathrm{d}t}{T} \quad \text{and} \quad \overline \nu \ \defeq \ \lim_{T \to \infty} \frac{\int_0^T \nu_T^\star \, \mathrm{d}t}{T}
\]
exist in $H^1(\Omega)$, since the optimal trajectories remain bounded and must satisfy
\[
\lim_{T \to \infty} \rho_T^\star(T) \ = \ \rho_\infty \quad \text{and} \quad \lim_{T \to \infty} \nu_T^\star(T) \ = \ \nu_\infty \; ,
\]
where the limits are understood in $H^1(\Omega)$. The latter claim follows by using an analogous argument as in the proof of Theorem~\ref{thm::ergodic}, recalling that the ergodic limit $(\rho_\infty,\nu_\infty)$ is strictly optimal on infinite horizons. This implies in particular that $0 = \mathcal A \overline \rho + \mathcal B \overline \nu$ as well as
\[
\mathcal L( \overline \rho, \overline \nu ) \ \leq \ \ell_\infty - \delta \ ,
\]
due to~\eqref{eq::contradiction} and Jensen's inequality, recalling that $\mathfrak L$ is convex. But this contradicts the definition of $\ell_\infty$ in~\eqref{eq::linf}. Thus, we must have
\[
\forall T > 0, \qquad \widetilde J(T,\cdot) \ = \ J_\infty^+(\cdot) \ = \ J_\infty(\cdot),
\]
where the latter equation follows from the fact the first equation holds uniformly and the fact that the contribution of the terminal cost,
\[
\lim_{T \to \infty} J_\infty^+(\rho_T^\star(T)) \ = \ 0
\]
vanishes for $T \to \infty$.

\bigskip
\noindent
\textit{Part III.} An immediate consequence of the constructions in Part II is that the unique weak solution of the parabolic HJB~\eqref{eq::HJBtilde},
\[
\forall x \in \Omega, \qquad \widetilde V(t,x) \ = \ V_\infty(x)                                                                                                                                                        \]
does not depend on $t > 0$. Thus, $V_\infty$ is a weak solution of the ergodic HJB~\eqref{eq::HJBinf} as claimed by the statement of this theorem. Finally, the energy estimate follows from~\eqref{eq::energy1} by using that
\begin{eqnarray}
\| V_\infty \|_{L^2(\rho_\infty)} & \ = \ & \sup_{\rho_0 \neq \rho_\infty} \frac{\int_{\Omega} V_\infty (\rho_0 - \rho_\infty ) \, \mathrm{d}x}{\| \rho_0 - \rho_\infty \|_{L^2(\rho_\infty^{-1})}} \notag \\[0.16cm]
&=& \sup_{\rho_0 \neq \rho_\infty} \frac{J_\infty^+(\rho_0)}{\| \rho_0 - \rho_\infty \|_{L^2(\rho_\infty^{-1})}} \notag \\[0.16cm]
&\overset{\eqref{eq::energy1}}{\leq}& \frac{1}{\gamma} \sqrt{\int_{\Omega} \ell[\mu_\infty]^2 \rho_\infty \, \mathrm{d}x} \; .
\end{eqnarray}
This completes the proof.\qed

\begin{remark}
The above framework can, in principle, be extended for systems with state constraints of the form
\[
\mathbb E\{ c(X(t)) \} \ \leq \ 0 \qquad \Longleftrightarrow \qquad \int_{\Omega} c \, \rho(t) \, \mathrm{d}x \ \leq \ 0
\]
for a constraint function $c \in L^\infty(\Omega)$. This can be done by adding the above linear inequality constraint to~\eqref{eq::ocp}, as long as one restricts the analysis to initial density functions $\rho_0$ for which this linear inequality constraint is strictly feasible, such that Slater's constraint qualification still holds. Because \eqref{eq::ocp} is then still a convex PDE-constrained optimization problem, it still admits a unique solution. One major difficulty, however, is that, in this case, $J(T,\cdot)$ is not a bounded linear operator anymore. This renders Lemma~\ref{lem::HJB} and Theorem~\ref{thm::HJB} inapplicable to systems with stochastic state constraints. Instead, the KKT optimality conditions of~\eqref{eq::ocp} have to be analyzed by using the Zowe-Kurcyusz theorem \cite{Zowe1979}, see also \cite[Thm.~1.56]{Hinze2009}. In combination with Rockafellar's theorem~\cite[Thm.~22.14]{Bauschke2011} it then follows that the Frechet subdifferential of $J(T,\cdot)$ (and also $J_\infty$) is, under suitable regularity assumptions, a maximal monotone operator. The practical relevance of such analysis techniques is, however, limited, as the whole point of the HJBs~\eqref{eq::HJB} and~\eqref{eq::HJBinf} is that one can compute a value function that is independent of $\rho_0$. Because this desirable property of the HJB fails to hold in the presence of explicit state constraints, a much more practical heuristic is to choose a function $\ell$ that takes large values in regions of the state space that one wishes to constraint. The optimal controller can then be expected to automatically avoid these regions with high probability whenever possible.
\end{remark}

\section{Global Optimal Control Algorithm}
\label{sec::algorithm}
The results of Lemma~\ref{lem::HJB} and Theorem~\ref{thm::HJB} motivate a completely new numerical method for solving HJBs. This method relies on discretizing the dynamic programming recursion for the bounded linear operator $J(T,\cdot)$ rather than on discretizing the value function $V$ directly.

\subsection{Galerkin Projection}
The bounded linear operator $J(T,\cdot)$, as defined in~\eqref{eq::ocp}, can be approximated by a finite dimensional linear function $\mathsf{J}_{0}(\cdot)$ by applying a simple Galerkin projection. This means that we choose basis functions
\[
\varphi_1, \varphi_2, \ldots, \varphi_\mathsf{m} \in H^1(\Omega)
\]
in order to discretize the states $\rho$  and controls $\nu$. For simplicity of presentation, we assume that this construction is based on given grid points
\[
\mathsf{x}_1,\mathsf{x}_2,\ldots,\mathsf{x}_{\mathsf{m}} \in \Omega \quad \text{and} \quad \mathsf{t}_k = kh
\]
with step-size $h > 0$ and time index $k \in \mathbb N$ and that Lagrange's equation holds\footnote{Note that the methods in this paper can be extended easily for all kinds of finite element-,  model order reduction, and other related numerical methods for linear parabolic PDEs or even tailored methods for FPKs, which might use basis functions that do not satisfy Lagrange's equation---this equation is merely introduced in order avoid unnecessary technical notation overhead.}; that is, $\varphi_i(\mathsf{x}_j) = \delta_{i,j}$, where $\delta_{i,j} = 0$ if $i \neq j$ and $\delta_{i,i} = 1$. In this case, the discrete approximation of the bounded linear operator $J(T,\cdot)$ at $T = t_\mathsf{N}$ has the form
\begin{eqnarray}
\label{eq::JN}
\mathsf{J}_{0}(z_0) \ \defeq \ &\min_{z,v} \ & \sum_{k=0}^{\mathsf{N}-1} \sum_{i=1}^{\mathsf{m}} \, \mathsf{w}_i \cdot \ell\left( \mathsf{x}_i, \frac{v_{k,i}}{z_{k,i}} \right) \cdot z_{k,i} \\[0.26cm]
&\text{s.t.} \ & \left\{
\begin{array}{l}
\forall k \in \{ 0,1, \ldots, \mathsf{N}-1 \},  \\[0.16cm]
\mathsf{E} z_{k+1} = \mathsf{A} z_k + \mathsf{B} v_k \\[0.16cm]
(z_{k,i},v_{k,i}) \in \mathbb K \quad \text{f.a.} \quad
i \in \{ 1, \ldots, \mathsf{m} \}.
\end{array}
\right. \notag
\end{eqnarray}
Here, $\mathsf{A}$ and $\mathsf{B}$ are discrete-time system matrices that are found by a Galerkin projection of $\mathcal A$ and $\mathcal B$ in combination with a suitable implicit differential equation solver in order to discretize the system also in time. In this context, $\mathsf{E}$ is an invertible mass matrix. The coefficients of the associated truncated series expansions
\[
\rho(\mathsf{t_k}) \ \approx \ \sum_{i=1}^\mathsf{m} z_{k,i} \varphi_i \quad \text{and} \quad \nu(\mathsf{t_k}) \ \approx \ \sum_{i=1}^\mathsf{m} v_{k,i} \varphi_i
\]
are optimized subject to a conic constraint, given by
\[
\mathbb K \ = \ \left\{ \
(z,v) \in \mathbb R \times \mathbb R^{n_u} \ \middle| \ v \in z \mathbb U \ \ \text{and} \ \ z > 0 \
\right\}.
\]
Moreover, the quadrature weights $\mathsf{w}_i > 0$ are needed in order to discretize the integral in the objective of~\eqref{eq::ocp}.

\begin{remark}
Methods for bounding the discretization error of the above Galerkin projection of~\eqref{eq::ocp} in dependence on $\mathsf{m}$ and $\mathsf{N}$ can be found in~\cite{Brenner2005}. Moreover, we recall that a desirable property of Galerkin methods is that---at least for appropriate choices of the basis functions---Markov's conservation laws
\[
\mathsf{1}^\tr \mathsf{E}^{-1} \mathsf{A} = \mathsf{1}^\tr \qquad \text{and} \qquad \mathsf{1}^\tr \mathsf{E}^{-1} \mathsf{B} = 0
\]
hold. Similarly, accurate discretization schemes either satisfy (or simply enforce) the condition
\[
\mathsf{A} + \mathsf{B} \cdot \mathrm{diag}( \mathsf{u} ) \ > \ 0
\]
componentwise for all vectors $\mathsf{u}$, with $\mathsf{u}_i \in \mathbb U$, such that $\mathsf{z}_{k} > 0$ implies $\mathsf{z}_{k+1} > 0$ for all feasible inputs.
\end{remark}

\subsection{Dual perspective function}
For the following developments it will be convenient to introduce the Fenchel conjugate of the stage cost,
\[
\mathsf{d}(x,\lambda) \ \defeq \ \min_{u \in \mathbb U} \, \left\{ \ell(x,u) + \lambda^\tr u \right\} \; .
\]
Fortunately, in practical applications, one can find an explicit expression for $\mathsf{d}$. For instance if $\ell$ is a diagonal quadratic form in $u$ and $\mathbb U$ a simple interval box, a single evaluation of $\mathsf{d}$ with modern processors can be done efficiently, usually, in a couple of nanoseconds. In this case, the dual perspective function
\[
\mathsf{D}(\lambda) \ \defeq \ \left(
\begin{array}{c}
\mathsf{d}( \mathsf{x}_1, \lambda_1/\mathsf{w}_1) \cdot \mathsf{w}_1 \\
\vdots \\
\mathsf{d}( \mathsf{x}_{\mathsf{m}}, \lambda_{\mathsf{m}}/\mathsf{w}_{\mathsf{m}}) \cdot \mathsf{w}_{\mathsf{m}}
\end{array}
\right)
\]
can be evaluated with reasonable effort, even if the number of basis functions, $\mathsf{m}$, is very large.

\subsection{Operator-Theoretic Dynamic Programming}
\label{sec::otdp}
Assuming that the stochastic optimal control problem is discretized as above, it can be solved by starting with $y_{\mathsf N} = 0$ and solving the recursion
\begin{align}
\label{eq::DP}
\mathsf{E}^\tr y_{k} \ = \ \mathsf{A}^\tr y_{k+1} + \mathsf{D}( \, \mathsf{B}^\tr y_{k+1} \, )
\end{align}
backwards in time, for $k \in \{ \mathsf N-1, \, \ldots, \, 1, \, 0 \}$. This procedure can be justified directly by observing that the linear cost-to-go function of the discrete optimal control problem~\eqref{eq::JN} is given by $\mathsf{J}_{k}(z) = y_k^\tr E z$. It satisfies the discretized operator-theoretic dynamic programming recursion,
\begin{eqnarray}
y_k^\tr \mathsf{E} z & \ = \ & \min_{v} \sum_{i=1}^{\mathsf{m}} \, \mathsf{w}_i \cdot \ell\left( \mathsf{x}_i, \frac{v_{i}}{z_{i}} \right) \cdot z_{i} + y_{k+1}^\tr ( \mathsf{A} z + \mathsf{B} v ) \notag \\
&=& \left[ \ \mathsf{A}^\tr y_{k+1} + \mathsf{D}\left( \, \mathsf{B}^\tr y_{k+1} \, \right) \ \right]^\tr z \; . \notag
\end{eqnarray}
A comparison of coefficients yields~\eqref{eq::DP}. Moreover, the associated separable minimizers
\[
u_{k,i}^\star \ \defeq \ \underset{u \in \mathbb U}{\mathrm{argmin}} \ \left\{ \ \ell(x_\mathrm{i},u) + \frac{y_{k+1}^\tr \mathsf{B}_i u}{\mathsf{w}_i} \ \right\}
\]
can be used to recover an approximation of the optimal control law $\mu^\star = \nu^\star/\rho^\star$ at the minimizer $(\rho^\star,\nu^\star)$ of~\eqref{eq::ocp},
\[
\mu^\star( \mathsf{t}_k, \mathsf{x}_i ) \ \approx \ u_{k,i}^\star \; .
\]
The computational complexity for computing this approximation of $\mu^\star$ is given by $\mathcal O( \mathsf{m} \cdot \mathsf{N} )$.
Finally, regarding the solution of the infinite horizon control problem, one option is to keep on iterating~\eqref{eq::DP} backward in time until a suitable termination criterion is satisfied, for example, until
\[
\| u_{k} - u_{k+1} \|_{\infty} \ \leq \ \mathsf{TOL}
\]
for a user specified tolerance $\mathrm{TOL} > 0$. This method converges for $k \to -\infty$ as long as the conditions of Theorem~\ref{thm::HJB} are satisfied---assuming that the above finite element discretization is reasonably accurate. Alternatively, one can use a large-scale convex optimization problem solver to solve~\eqref{eq::JN} for $\mathsf{N}=1$ but subject to the additional stationarity constraint $z_1 = z_0$. 

\begin{remark}
As reviewed in the introduction, the duality between HJBs and optimally controlled FPKs is known since a long time and, in fact, in~\cite{Rutquist2017} a method for solving HJBs via FPK-constrained convex optimization can be found. The operator-theoretic dynamic programming recursion~\eqref{eq::DP} is, however, an original contribution of this paper. It is entirely explicit and does---at least if the Fenchel conjugate $d$ admits an explicit expression---not even require an optimization problem solver at all. 
\end{remark}

\begin{figure*}[t]
\centering
\includegraphics[scale=0.15]{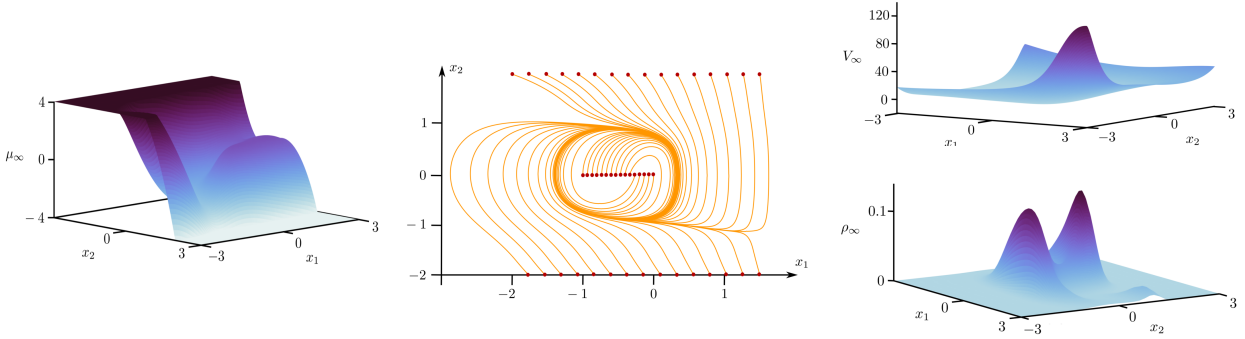}
\caption{\label{fig::illustration} \texttt{LEFT:} the optimal ergodic control law $\mu_\infty: \Omega \to \mathbb U$ on the domain $\Omega = [-3,3]^2$ for the control bounds $\mathbb U = [-4,4]$. \texttt{MIDDLE:} simulations of the states of the closed-loop system without noise for selected initial values (red dots), which converge to a non-trivial attractor set. \texttt{UPPER RIGHT:} the indefinite solution $V_\infty$ of the ergodic HJB that corresponds to the unique Riesz representation of the cost function, as introduced in Theorem~\ref{thm::HJB}. \texttt{LOWER RIGHT:} the optimal multi-modal steady state $\rho_\infty$.}
\end{figure*}

\subsection{Numerical Illustration}
This section is about a nonlinear stochastic system with two states and one control input, given by
\[
f(x) = \left(
\begin{array}{c}
x_2 - \frac{1}{2} x_1 x_2 \\
-x_1
\end{array}
\right) \quad G(x) = \left(
\begin{array}{c}
0 \\
1
\end{array}
\right), \ \ \epsilon = 0.2,
\]
on the domain $\Omega = [-3,3]^2$ and control constraint set $\mathbb U = [-4,4]$. An associated non-convex stage cost function is given by
\[
\ell(x,u) = \frac{1}{4}x_1^2 + 3 \left( x_2^2 - 1 \right)^2 + \frac{1}{2} u^2 \; .
\]
This example trivially satisfies Assumptions~\ref{ass::blanket},~\ref{ass::Q}, and~\ref{ass::VAR}. Assumption~\ref{ass::LMI} can be verified numerically finding the exponential decay rate $\gamma \approx 0.06$. Figure~\ref{fig::illustration} visualizes the optimal infinite horizon control law $\mu_\infty$, the optimal steady state $\rho_\infty$, and the associated dual solution $V_\infty$ that solves the ergodic HJB at
\[
\ell_\infty \ \approx \ 2.03 \; .
\]
This solution has been found by using a Galerkin discretization with Lagrange basis functions on a simple equidistant grid with $\mathsf{m} = 22500$ basis functions and $N = 10^3$ time iterations with $h = 0.05$.  The table below,
\begin{center}
\begin{tabular}{|c|c|c|c|}
\hline
$\mathsf{m}$ & $\mathsf{N}$ & $\| \mu_\infty^{\mathsf{m}} - \mu_\infty^\star \|_{L^2(\Omega)}$ & $\mathsf{CPU-time}$ \\
\hline
$100$ &  $10^3$ & $0.17$ & $9 \, \mathrm{ms}$ \\
\hline
$400$ &  $10^3$ & $0.015$ & $35 \, \mathrm{ms}$ \\
\hline
$900$ &  $10^3$ & $0.0056$ & $79 \, \mathrm{ms}$ \\
\hline
\end{tabular} \ ,
\end{center}
lists the CPU run-time of a \texttt{JULIA}-prototype implementation of~\eqref{eq::DP} for different choices of $\mathsf{m}$. The third column of this table lists the $L^2$-norm between the control law $\mu_\infty^{\mathsf{m}}$ that is found by using $\mathsf{m}$ basis functions and the control law $\mu_\infty^\star$ that is found by using $22500$ basis functions. 

\bigskip
\noindent
Note that this example should be regarded as a ``toy problem''. More advanced implementations of the above method and case studies are left for future work. Nevertheless, the above example is non-trivial in the sense that it has a multi-modal steady state $\rho_\infty$. Moreover, the associated solution $V_\infty$ of the ergodic HJB turns out to be indefinite. And, finally, the closed-loop system has a non-trivial limit behavior: the plot in the center of Figure~\ref{fig::illustration} shows closed-loop simulations without noise in order to illustrate that the nominal trajectories converge to a non-trivial periodic attractor.

\section{Conclusions}
\label{sec::conclusions}
This paper has presented methods for analyzing ergodic- and infinite horizon stochastic optimal control problems. In detail, Theorem~\ref{thm::ergodic} has established the existence and uniqueness of ergodic optimal controls under Assumptions~\ref{ass::blanket} and~\ref{ass::Q}. These assumptions are satisfied for most nonlinear stochastic control problems of practical relevance, because nonlinear McKean-Vlasov systems on bounded domains with reflecting boundaries always admit a Lyapunov-Has'minski\v{i} function. But even if the domain is unbounded such functions \mbox{are---in} contrast to traditional Lyapunov functions---relatively easy to find.

\bigskip
\noindent
Concerning the analysis of infinite horizon optimal control problems, this paper has proposed to introduce a novel class of Lyapunov diffusion operators that can be used to formulate a generalized Bakry-Emery condition, as summarized in Assumption~\ref{ass::LMI}. The relevance of this Lyapunov diffusion operator based LMI has been clarified in Theorem~\ref{thm::HJB}, which has used this condition to establish the existence of solutions to ergodic HJB equations for nonlinear optimal control problems with indefinite stage costs.

\bigskip
\noindent
For the first time, this paper has presented methods for analyzing a very general class of nonlinear HJBs by using functional analysis methods in Hilbert-Sobolev spaces that would usually find applications in the context of linear parabolic PDEs. This has far reaching consequences on the development of global robust optimal control algorithms. This is because the proposed operator-theoretic dynamic programming method renders existing numerical methods for discretizing linear parabolic PDEs fully applicable for solving nonlinear HJBs. The exploitation of the implications of this observation appear to be a fruitful direction for future research on stochastic optimal control.

\bibliographystyle{plain}
\bibliography{references}

\end{document}